\documentclass[3p,11pt,authoryear]{elsarticle}

\usepackage{amssymb,graphicx,hyperref}
\usepackage{mathtools}
\usepackage{appendix}
\usepackage{color}
\usepackage{subfigure}
\usepackage{bm}
\usepackage{dsfont}
\usepackage{tabularx}
\usepackage{ltablex}
\usepackage{amsmath}
\usepackage{mathrsfs,BOONDOX-cal}
\usepackage{enumitem}
\usepackage{cases}
\usepackage{empheq}
\usepackage{diagbox}
\usepackage{arydshln}
\usepackage{soul}

\usepackage[linesnumbered,lined,boxed]{algorithm2e}
\usepackage{setspace}
\SetAlCapSkip{0.5em}
\SetKwComment{Comment}{$\triangleright$\ }{}

\newcommand{\tabincell}[2]{\begin{tabular}{@{}#1@{}}#2\end{tabular}}

\newtheorem{proposition}{\bf Proposition}

\newtheorem{lemma}{\bf Lemma}

\newcounter{example}
\newenvironment{example}[1][]{\refstepcounter{example}\par\medskip
   \noindent 
{\bf Example~\theexample:}~#1\it}{\medskip}

\newcommand{\proof}[1]{{\it #1}}

\newcommand{\Halmos}{{\hfill$\blacksquare$}}
\usepackage{natbib}
 \bibpunct[, ]{(}{)}{,}{a}{}{,}%

\def\map{\Gamma} 
\def\mL{\mathcal{L}} 
\def\mE{\mathcal{E}} 

\graphicspath{{./figs/}}

\journal{European Journal of Operational Research}

\begin{document}

\begin{frontmatter}

\title{A Restless Bandit Model for Dynamic Ride Matching with Reneging Travelers}

\author[jing]{Jing Fu}
\affiliation[jing]{organization={School of Engineering, RMIT University},
            city={Melbourne},
            postcode={VIC3000}, 
            state={Victoria},
            country={Australia}}
    
\ead{jing.fu@rmit.edu.au} 

\author[joyce]{Lele Zhang}
\affiliation[joyce]{organization={School of Mathematics and Statistics, The University of Melbourne},
            city={Melbourne},
            postcode={VIC3010}, 
            state={Victoria},
            country={Australia}}

\author[zhiyuan]{Zhiyuan Liu}
\affiliation[zhiyuan]{organization={School of Transportation, Southeast University},
            city={Nanjing},
            state={Jiangsu},
            country={China}}  

\begin{abstract}
This paper studies a large-scale ride-matching problem with a large number of travelers who are either drivers with vehicles or riders looking for sharing vehicles. 
Drivers can match riders that have similar itineraries and share the same vehicle; and reneging travelers, who become impatient and leave the service system after waiting a long time for shared rides, are considered in our model.
The aim is to maximize the long-run average revenue of the ride service vendor, which is defined as the difference between the long-run average reward earned by providing ride services and the long-run average penalty incurred by reneging travelers.
The problem is complicated by its scale, the heterogeneity of travelers (in terms of origins, destinations, and travel preferences), and the reneging behaviors.
To this end, we formulate the ride-matching problem as a specific Markov decision process and propose a scalable ride-matching policy, referred to as Bivariate Index (BI) policy.
The BI policy prioritizes travelers according to a ranking of their bivariate indices, which we prove, in a special case, leads to an optimal policy to the relaxed version of the ride-matching problem. 
For the general case, through extensive numerical simulations for systems with real-world travel demands, it is demonstrated that the BI policy significantly outperforms baseline policies. 
\end{abstract}


\end{frontmatter}
\section{Introduction}\label{sec:introduction}

Traffic congestion in large cities has caused significant problems in the environment, human society and economy. In Australia, congestion cost has been expected to be over 27.7 billion Australian dollars by 2030 (see \citep{AAA}). Meanwhile, the average occupancy of private vehicles (i.e., the number of people per car) is very low, around 1.2 people per car in Melbourne (see \citep{vicroads}) and slightly above 1.2 across Australia. 
By improving vehicle occupancy, proactive techniques have a great potential to reduce (the growth of) traffic demands and mitigate congestion.

Ride-sharing services, offering shared rides for travelers who have similar itineraries, are an effective proactive technique to improve the utilization of vacant seats and increase the vehicle occupancy. From \citep{dOrey2012}, the average occupancy per travel kilometer would increase by 48\% once providing ride-sharing services. 
Ride-sharing enables fare-sharing between travelers with door-to-door services and comfortable riding experiences: it balances the traveling experiences and financial benefits of travelers. 
Also, \cite{Wang2018} found that ride-sharing can reduce 20\%-30\% travel distance and 20\%-25\% taxi requirement in Singapore.

Methodologies for ride-sharing services have been considered in accordance with different respects. 
For instance, for systems with determined riders for each vehicle, conventional routing methods with enabled time window management have been developed in \citep{Cordeau2003, Baldacci2004}. Dynamic posted pricing mechanisms, such as auction-based and negotiation-based mechanisms, were also used in \citep{Kleiner2011, Furuhata2013, Nourinejad2016} 
 to improve the quality of ride-sharing services, for which the itinerary prices usually varied with the number of travelers, and travel distance and trip urgency. 
The service quality can be further improved by appropriately \emph{matching} different travelers with similar itineraries, as studied in \citep{Ma2014, Wang2018}, which is the main concern of this paper.

The rider-matching modeling and optimization have been widely studied for years. 
\cite{Hong2017} and \cite{Dong2018} considered ride-matching optimizations in a \emph{static} manner, where traveler demands were required to be provided a long time in advance. 
\cite{guo2021real} and \cite{guo2023modelling} studied ride-sharing optimizations with given sets of riders and drivers for each time slot and aimed at the local optimum in that time slot.
Rapid developments of emerging IT techniques have also enabled real-time ride-sharing services for mobile applications with relatively long-term objectives. It followed a large number of studies involving \emph{dynamic} optimizations, which typically re-considered the entire system decisions upon each newly arrived ride demand and suffered high computational complexities as system scales became large.
In particular, \cite{Lokhandwala2018} studied the New York City taxi network using an agent-based model that incorporated different traveler preferences. They demonstrated that the vehicle occupancy can be significantly improved by the ride-sharing mechanism.  
\cite{Hosni2014} proposed a mixed integer programming model to dynamically arrange new passengers to share rides with on-road people.
\cite{Wang2018, Ma2014} developed matching strategies and considered specific features of ride-matching applications, such as pick-up and drop-off windows, passenger level agreements and preferences, and delayed arrivals. 
We focus on a ride-matching problem that dynamically optimizes newly arrived travelers without requesting the actual demand a long time a priori.

More importantly, travelers/passengers are likely to be impatient while waiting for someone who can match their rides: the travelers are usually allowed to change their minds and choose other (non-sharing) alternatives (see \citep{Shi2016}).
We refer to these impatient travelers as the \emph{reneging} travelers as in \citep{Pazgal2008}.
The reneging phenomenon plays an important role in real-world applications (see \cite{Kim2016}), and to the best of our knowledge, most existing studies, such as \citep{Mamalis2019}, assume travelers are truthful and do not renege. 
In this paper, we take into account impatient travelers; if they go for other services while waiting for a shared ride, a certain penalty will be generated as negative feedback for the service vendor.

In \citep{Conolly2002, Afeche2014}, reneging behaviors were studied in two-sided queueing models with exponentially distributed inter-reneging time, where the strategies for matching travelers (over multiple queues) and routing taxis were re-optimized when a traveler arrived or reneged.
It was found that the two-sided queueing model can address the taxi service, since taxi drivers do not have designated destinations or preferences in choosing travelers.
\cite{agussurja2019state} studied a last-mile ride-sharing problem through a discrete-time Markov decision process (MDP) and assumed a two-time-slot waiting period for each passenger. If a passenger waits for over two time slots, he or she is considered to renege. 
In this paper, we do not assume any maximal waiting time of riders or drivers and allow both the riders and drivers to have different origins (i.e., sources), destinations, and ride preferences, which complicates the entire problem and prevents the above-mentioned methods from being applied directly.

Appropriate reneging control can further help reduce the average waiting time of travelers, which has been widely considered as an important criterion for service quality in previous studies such as \citep{Wardman2001}.
Here, the traveler waiting time is restricted and reflected by imposing the penalties incurred by reneging passengers as negative rewards attached to our system. 
By maximizing the long-run average total reward, which accounts for the reneging penalties, we aim to achieve a trade-off between the monetary revenue gained through the shared rides and traveler experiences.
A detailed definition of our objective is provided in Section~\ref{sec:model}.

Our problem is significantly complicated by the heterogeneity of travelers, their reneging behaviors and the large scale of systems.
We model it through as an MDP with high-dimensional state space and approximate its optimality by applying the Whittle relaxation technique (see \citep{whittle1988restless}), which was originally proposed for the restless multi-armed bandit (RMAB) problem.

The ride-matching problems have been studied through MDP-based methodologies in recent years.
\cite{xu2018large} adopted the Temporal-Difference (TD) learning algorithm for dynamic pricing in the ride-hailing problem, where the locations, profiles and environment conditions of the drivers were modeled through MDPs. 
The TD-based method for ride-hailing was further considered by \cite{chen2019inbede,li2019efficient}.
Later, \cite{tang2019deep} studied a similar ride-hailing pricing problem, for which the value functions of the underlying MDPs were learnt through the deep neural network.
In the work by \cite{xu2018large,chen2019inbede,li2019efficient,tang2019deep}, ride-matching was also considered through maximizing the instantaneous utility for the learnt value functions in each time slot, where order cancellation or rider/passenger reneging was failed to be formulated.
\cite{qin2021optimizing} and \cite{ke2022learning}
optimized ride-matching problems, taking into account the waiting and reneging of riders, by deciding whether or not each driver should take an order based on a given matching strategy.
The given matching strategy, similar to the previous work, was obtained through optimizing an instantaneous objective in each time slot.
\cite{ozkan2020dynamic} formulated a ride-sharing problem in the manner of a bipartite tree with arriving and reneging riders and drivers.
\cite{ozkan2020dynamic} focused on the maximization of the long-run sum of the weighted match rates by appropriately selecting drivers for the arrived riders.
Unlike \cite{ozkan2020dynamic}, 
in this paper, we consider selections for both drivers and riders/passengers, leading to possible reservation of some drivers for not-yet-arrived riders who are more profitable than those currently waiting in the system.

\cite{whittle1988restless} proposed the Whittle relaxation technique to analyze the RMAB problem, which exhibits high-dimensional state space and
has been proved to be PSPACE-hard in general (see \citep{papadimitriou1999complexity}). 
A standard RMAB problem consists of parallel binary-action MDPs, referred to as bandit processes, that are coupled by constraints over their action variables.
Each bandit process is dependent on a single constraint for the conventional RMAB problem, such as the problems discussed in \cite{whittle1988restless,nino2007dynamic,graczova2014generalized} and \cite{fu2021large}, and the extended RMAB problems with multi-actions studied in \cite{glazebrook2011general,hodge2015asymptotic,zayascaban2019asymptotically,killian2021qlearning} and \cite{xiong2022reinforcement}; while, here, the ride-matching problem involves multiple action variables for each bandit process that are subject to multiple action constraints.
It follows with stronger dependencies between the bandit processes and prevents conventional techniques from quantifying the marginal rewards/costs for matching different travellers, which is a key step for analyzing RMAB or some extended versions of RMAB problems.
A detailed explanation of the difference between the ride-matching and the past work about the RMAB problems is provided in Section~\ref{subsec:graph}.

A preliminary conference version of this paper was presented
in \citep{fu2021large}. Here, we extend the conference work by providing a more thorough understanding of the proposed heuristic policy as well as an analytical analysis of its properties, taking into consideration of multiple travel origins, drivers' eligibility and time-inhomogeneous arrival rates of travelers. 

\vspace{-0.2cm}
\subsection*{Main Contributions}

We formulate the ride-matching problem with hundreds of thousands of daily travelers, who have potentially different travel sources, destinations, preferences and means, as an MDP consisting of parallel stochastic processes.
The stochastic processes are again MDPs, referred to as \emph{sub-processes}, and are correlated with each other through multiple constraints.
We adopt and extend the Whittle relaxation technique and propose a policy for the ride-matching problem, 
for which the storage and computational complexity of implementing the proposed policy is linear and logarithmic, respectively, to the number of candidate matches between different riders and drivers.

As mentioned above, compared to the past work on RMAB or extended versions of RMAB problems, the ride-matching problem exhibits stronger dependencies between the sub-processes, preventing past methodologies from being directly applied.
We propose a method that summarizes the variate of the Lagrange multipliers associated with the multiple action constraints applied to the same sub-process.
Through this method, we quantify the marginal profits of selecting certain matches between riders and drivers by real numbers, which are referred to as the \emph{bivariate indices} of the candidate matches.
In the proposed policy, we prioritize travelers according to their \emph{bivariate indices} which evolve dynamically due to the time-varying system states.

We prove that, {\color{black}when the reneging rates are negligible and arrival rates of the drivers and riders are compatible to each other}, the bivariate indices lead to a policy optimal to the relaxed version of the ride-matching problem.
For the general case, the effectiveness of the proposed policy is demonstrated through extensive numerical simulations with real-world settings: compared to the baseline policies, it achieves more than $25\%$ higher average rewards with compatible average traveler waiting time.

To the best of our knowledge, no existing techniques that considered traveling preferences of both drivers and riders and their reneging behaviors can be directly applied here; nor there is existing work that discussed the sub-process coupled by multiple action constraints.

The rest of the paper is organized as follows. Section~\ref{sec:model} describes the ride-matching problem formulated as an MDP with a stochastic optimization model. Section~\ref{sec:policies} proposes a bivariate index policy, an extension of the Whittle index policy, to solve the optimization problem. Numerical tests are conducted on a generic network with uniform service demand and a real-world network with heterogeneous demand to evaluate the performance of the proposed index policy in Section~\ref{sec:example}. Section~\ref{sec:conclusion} concludes the paper with a brief discussion.

\vspace{-0.3cm}
\section{A Ride-Matching Problem}\label{sec:model}
\vspace{-0.2cm}

We use $\mathbb{N}_{+}$ and $\mathbb{N}_{0}$ as the sets of positive and non-negative integers respectively.
Let $\mathbb{R}$, $\mathbb{R}_{+}$ and $\mathbb{R}_{0}$ be the sets of all, positive and non-negative reals, respectively.
We use $\backslash$ for the exclude operation for sets.

\vspace{-0.2cm}
\subsection{Ride-Sharing System}\label{subsec:ride-sharing}

Consider a ride-sharing system with $L$ different types of travelers geographically deployed and awaiting transport services. 
Here, a traveler is a logic term that represents either one person or a group/family of people.
Among these travelers, $L^D$ types of travelers have their own vehicles while the other $L^R \coloneqq L-L^D$ types of travelers have no vehicle and can only share a ride with the travelers with vehicles.
We refer to them as the \emph{drivers} and \emph{riders}, respectively. 
{\color{black}Let $\mL$ denote the set of all types of travelers and $\mL^D$ and $\mL^R$ denote the sets of drivers and riders respectively. Then $L=|\mL|$ is the cardinality of the set $\mL$, and analogously we let $L^D=|\mL^D|$ and $L^R=|\mL^R|$.}
We refer to a traveler of type $\ell\in{\color{black}\mL}$ as an $\ell$-traveler.
Drivers can share their vehicles with riders for some monetary profits, while riders wish to pay for a ride with the drivers.
Both drivers and riders may leave the system when they become impatient.
Travelers of the same type will originate and destine to the same geographical regions, such as a neighborhood or a suburb, in which all addresses are reachable from each other.
Travelers of different types differ from others in terms of origins, destinations, arrival rates, possession of vehicles, \textit{reneging rates} and potentially other travel preferences.

The travelers of type $\ell\in{\color{black}\mL}$ keep arriving in a Poisson process with rate $\lambda_{\ell}$, which is referred to as the arrival rate of type $\ell$. 
If, for a given period of time, a type of travelers arrive independently and identically and the number of arrivals is sufficiently large for the interested stochastic process to become stationary, then it is reasonable to model the arrival process with a time-homogeneous Poisson process.
In particular, for the more general case, we consider a time-inhomogeneous Poisson process consisting of successive exogenous time intervals where the Poisson rate varies from one interval to the next (see \citep{menon2017predicting, kieu2018stochastic}). 
The time-inhomogeneous Poisson process is able to capture busy and idle periods of transport networks and has been widely used for modeling traveler demands in the literature (see \cite{anwar2013changinow,matias2019method}).
Our analysis of system performance with arrival rates $\lambda_{\ell}$ for a given time interval can be applied accordingly to all successive intervals with highly diverse traveler demands.
In Section~\ref{sec:example}, we will apply the analytical results and numerically demonstrate the system performance in practical situations with time-varying arrival rates.

A traveler, either driver or rider, may become impatient and leave the system after waiting for a long time. The reneging rate describes the reciprocal of the expected waiting time after which 
travelers become impatient and depart the ride-sharing system without getting served.
Travelers with higher reneging rates will become impatient more easily and are more sensitive to waiting time.
We consider a case where the waiting time, for which a traveler of type $\ell_1\in{\color{black}\mL}$ becomes annoyed when awaiting a shared ride with another traveler of type $\ell_2\in{\color{black}\mL}$, is exponentially distributed with
a rate dependent on both $\ell_1$ and $\ell_2$. We refer to this rate as the reneging rate.
We postpone the detailed definitions of the reneging rates for a pair of matched travelers and start with the exposition of the matching mechanism of the travelers.

A traveler arrives in the system with a specific origin and destination and seeks a shared ride leading there.
We assume that two such travelers can and are agreed to share one vehicle for traveling if one is a driver and the other is a rider and their traveling routes to some extent \textit{match} each other.

In this paper, we focus on scheduling policies that match riders to drivers that are eligible to take a shared ride.
In Sections~\ref{sec:model} and \ref{sec:policies}, for the general case, the eligibility for a driver to match a rider is decided by a given rule.
The rule can be arbitrary if, for any eligible driver-rider pair, there exists a route from the origin of the driver to the origin of the rider and a route between the destination of the rider to the destination of the driver. 
Different eligibility rules may lead to different efficiency and performance of the same scheduling policy, which is about future work and is out of the discussion of this paper.
Proposition \ref{prop:priority} is applicable to any such rule for the match eligibility.
In Section~\ref{sec:example}, for the numerical simulations, we consider specific settings, for which the match eligibility is decided by a rule related to the travel distances of the shared rides and the gained revenue of the service vendor.
Here, we provide in Example~\ref{example:1} a simple example of the matching mechanism and a scheduling policy.

\begin{figure}[t]
\centering
\subfigure[]{\includegraphics[width=0.35\linewidth]{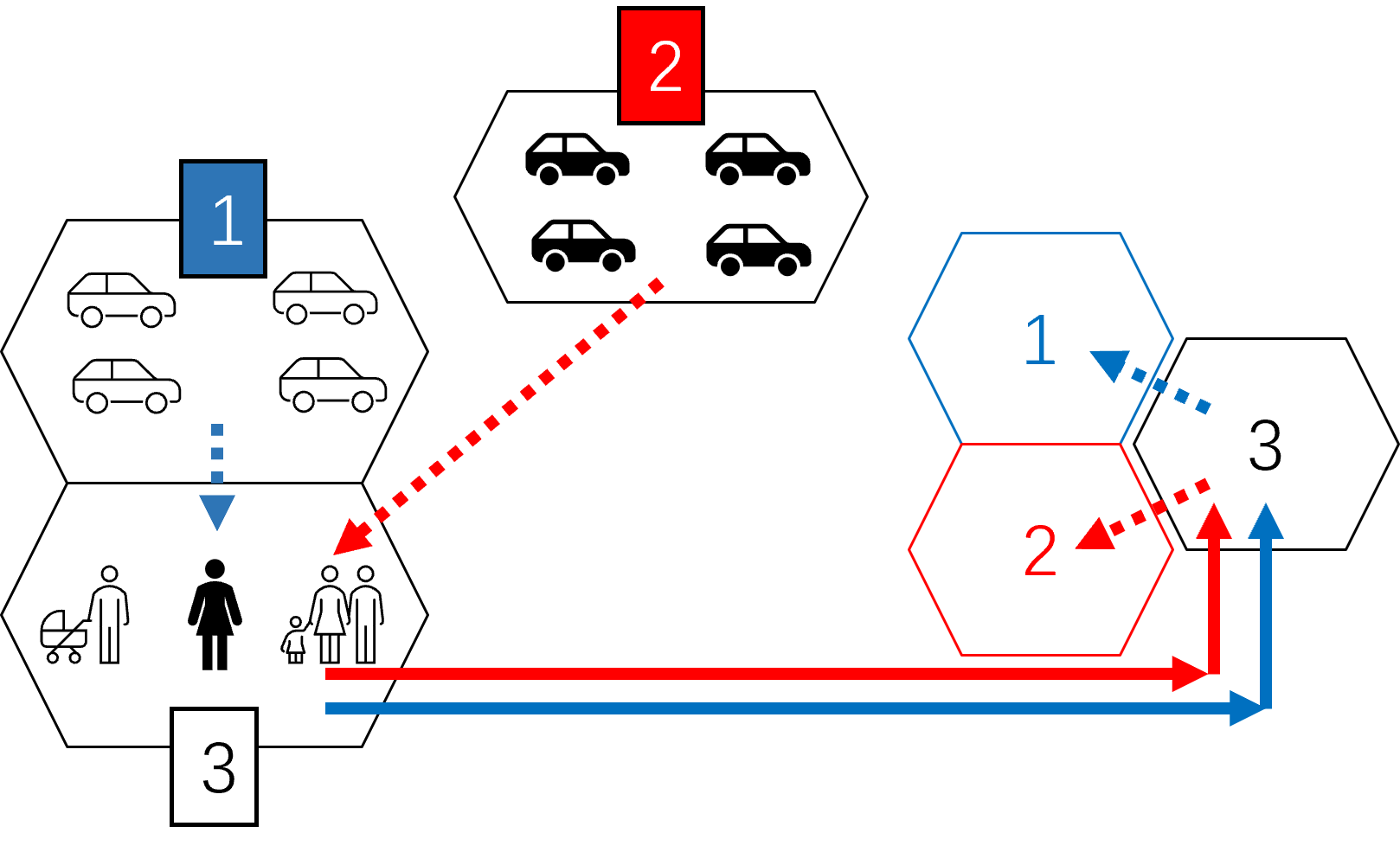}\label{fig:example1}}
\subfigure[]{\includegraphics[width=0.35\linewidth]{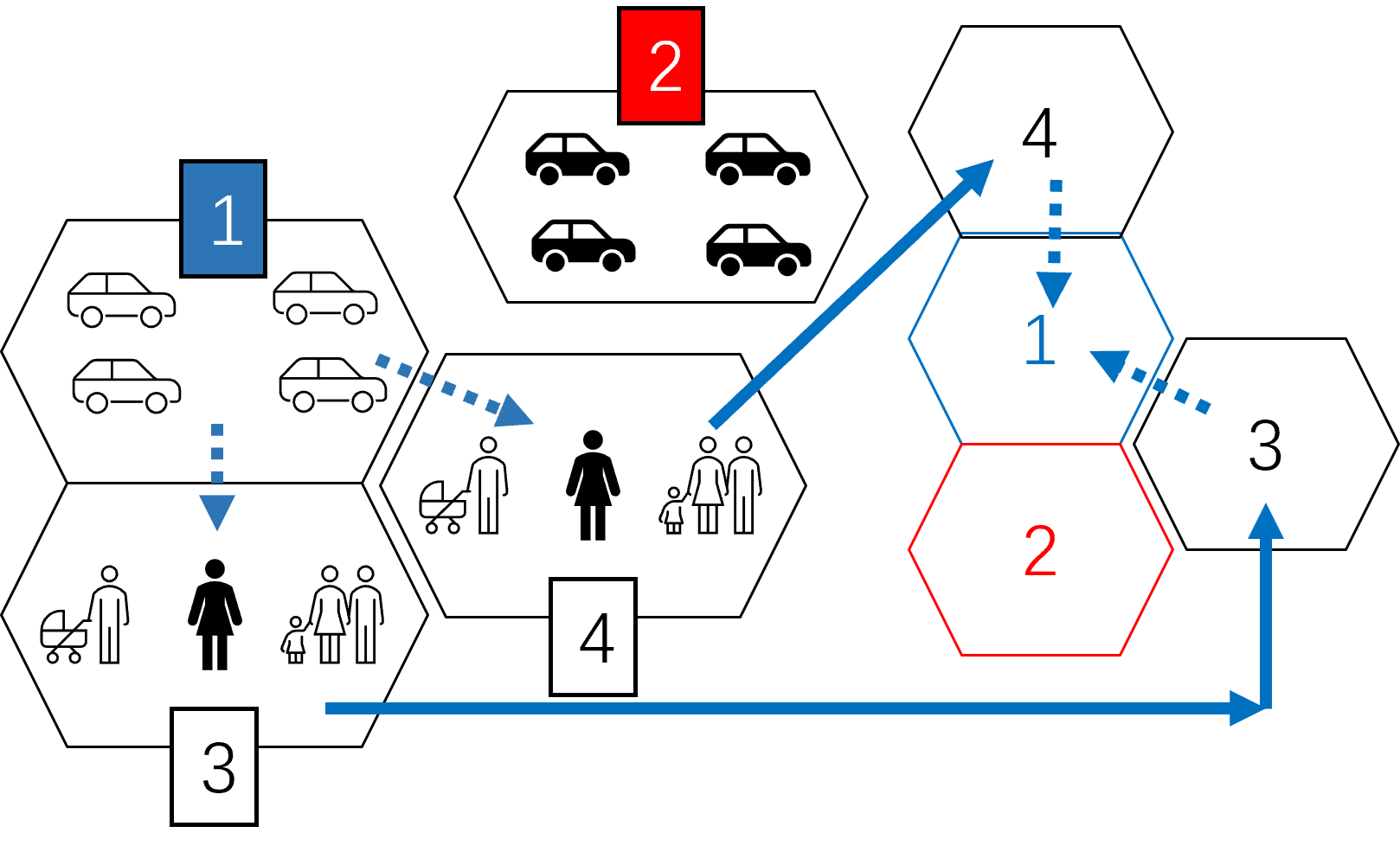}\label{fig:example2}}\vspace{-0.5cm}
\caption{A simple example for the ride-sharing system and the scheduling policies.}
\end{figure}

\begin{example}\label{example:1}
In Figure~\ref{fig:example1}, consider $L=3$ traveller types, including $L^D = 2$ types of drivers and one type of riders. 
That is, travellers of types 1 and 2 are drivers, and type-3 travellers are riders.
The hexagons with white and black cars represent the origin areas of type-1 and type-2, respectively, riders. The hexagon with some people is the origin area of the riders.
The hexagons labeled 1-3 on the right part of Figure~\ref{fig:example1} are the destination areas of travellers of types 1-3, respectively.
In this example, based on a given rule for match eligibility, 
both types of drivers are eligible to match the riders. 
Upon the arrival of a rider (driver), an employed scheduling policy will select a type of drivers (riders) to match the rider (driver) or, if the total number of waiting riders (drivers) exceeds a given threshold, directly reject the rider (driver) for ride sharing. 
The given threshold can be an arbitrarily large integer, which is used to avoid infinitely long waiting queues and to keep the entire system stable. 
There is a rigorous definition of the scheduling policy in Section~\ref{subsec:graph}.
In this example, if a type-1 driver is decided to match a rider, then the driver will go from its origin area to the origin area of the rider, as indicated by the dashed blue arrow in Figure~\ref{fig:example1}.
After the rider is picked up at the origin, the driver will firstly take the rider to the rider's destination, illustrated by the solid blue arrow in Figure~\ref{fig:example1}, and then go to its own destination as illustrated in Figure~\ref{fig:example1} by the dashed blue arrow from hexagon 3 to hexagon 1.  
Similarly, if a type-2 driver matches a rider, then the driver/rider will go as indicated by the red arrows in Figure~\ref{fig:example1}.

Different riders of the same type can be potentially matched by different types of drivers due to the dynamics of the system.
For instance, in this example, upon the arrival of a rider, if there is no type-1 driver but a lot of type-2 drivers waiting in the system, it may be sensible to match the rider with a type-2 driver. 
However, 5 minutes later, all the type-2 drivers have taken some previous riders.
The employed scheduling policy may think it is time to tell the newly arrived rider to wait for a type-1 driver.

On the other hand, each driver can potentially be eligible to take more than one type of riders. 
In Figure~\ref{fig:example2}, we add one more type of riders -- the type-4 riders. In Figure~\ref{fig:example2}, these riders originate from the hexagon with people and label 4 and destine to the right hexagon labelled 4.
Based on the given rule of match eligibility, type-1 drivers are eligible to take riders of both types 3 and 4.
In this context, upon the arrival of a type-3/4 rider, the employed scheduling policy needs to take decisions about whether a type-1 driver will take the shared ride.
\end{example}

When a driver matches a rider, we refer to the journey from the driver leaving its origin, going through the origin and the destination of the rider, until the driver reaches its destination as a \emph{shared ride}.
Upon the arrival of a type-$\ell$ rider, if the scheduling policy matches it with type-$\ell'$ drivers, and there is coincidentally a type-$\ell'$ driver awaiting type-$\ell$ riders, then the driver leaves its origin and starts the shared ride immediately. 
If, unfortunately, there is no type-$\ell'$ driver awaiting type-$\ell$ riders, then this newly arrived type-$\ell$ rider will wait in the system until a new type-$\ell'$ driver comes and is matched to it or reneges. 
We allow scheduling policies that match a rider to a not-yet-arrived driver, although such policies may be sub-optimal.
Upon the arrival of a driver, the employed scheduling policy will make decisions about how to match it with a rider.
Similarly, the driver will start the shared ride immediately if there is a matched rider awaiting it; otherwise, it waits in the system until a new rider arrives and matches it or reneges.

For each shared ride, the driver will go to the origin of the rider, the destination of the rider and its own destination, successively.
Such a travelling route can be generated in various manners; for instance, it can either be the shortest route to the rider's destination continued with the shortest route to the driver's destination, or be a combination of the shortest route and the longest route to the two destinations, respectively. 
Although the latter may not be a rational choice. 
We potentially allow travelers to specify their preference of traveling routes, providing the revenues gained by service vendors are always positive.

For each shared ride, the driver and/or the rider may change their minds with certain probabilities before reaching their destinations. For example, the rider may drop mid-way. 
For such an unsuccessful ride, the driver and/or the rider may join the ride-sharing system again with different origins, destinations, travel preferences et cetera. 
In this paper, the re-joined drivers and/or riders are modeled through the arrival processes of all the traveler types mentioned in earlier paragraphs and will be considered as new arrivals of appropriate traveler types.
In Sections~\ref{sec:model} and \ref{sec:policies}, the theoretical discussions and the proposed policies are directly applicable to the scenarios with positive probabilities of unsuccessful rides.
The probabilities of unsuccessful rides for different shared rides will affect the expected profits gained by selecting the driver-rider matches associated with these rides, which will be reflected in the long-run objective of the optimization problem.

\subsection{Stochastic Model}\label{subsec:graph}

Define an undirected bipartite graph $G = ({\color{black}\mL^D}, {\color{black}\mL}\backslash{\color{black}\mL^D},\mE)$ with edges between drivers and riders, where ${\color{black}\mL^D}$ is the set of all driver types,  ${\color{black}\mL}\backslash{\color{black}\mL^D}$ is the set of all rider types, and $\mE$ is the set of the edges. If there is an edge $(\ell_1,\ell_2)\in \mE$ with $\ell_1\in{\color{black}\mL^D}$ and $\ell_2\in{\color{black}\mL}\backslash{\color{black}\mL^D}$, then drivers of type $\ell_1$ and riders of type $\ell_2$ are eligible to match each other.
We refer to an edge $(\ell_1,\ell_2) \in \mE$ as an eligible \emph{match} between a driver of type $\ell_1$ and a rider of type $\ell_2$.
We label all the edges in $\mE$ by integers $e \in {\color{black}\mE}$, where $E=|\mE|$ is the cardinality of the set $\mE$, alternatively refer to an edge labeled by $e\in {\color{black}\mE}$ as the edge $e$. For any edge $e\in{\color{black}\mE}$, let $\ell^D(e)$ and $\ell^R(e)$ represent the vertices incident to the edge $e$ in sets ${\color{black}\mL^D}$ and ${\color{black}\mL}\backslash {\color{black}\mL^D}$, respectively.
For any $\ell\in{\color{black}\mL}$, define $\mE_{\ell}\coloneqq \{e\in{\color{black}\mE}| \ell^D(e)= \ell \text{ or } \ell^R(e)=\ell\}$ as the set of all the eligible matches for traveler type $\ell$ (all the edges incident to the vertex $\ell$). Assume that, for all $\ell\in{\color{black}\mL}$, $\mE_{\ell}\neq \emptyset$.

The travelers, including both drivers and riders, will keep arriving in the ride-matching system, and the system controller will assign a traveler of type $\ell\in{\color{black}\mL}$, once they arrive, to be \emph{served} by a match in $\mE_{\ell}$. 
For $\ell\in{\color{black}\mL}$ and $e\in {\color{black}\mE}$, let $N_{\ell,e}(t)\in\mathbb{N}_0$ represent the number of travelers of type $\ell$ that are waiting to be matched by travelers of another type, of which the corresponding vertex is incident to the edge $e$ ($e\in\mE_{\ell}$), at time $t\geq 0$. Define, without loss of generality, $N_{\ell,e}(t)\equiv 0$ if $e\notin \mE_{\ell}$, and assume an initial condition $N_{\ell,e}(0) = 0$ for all $\ell\in{\color{black}\mL}$ and $e\in{\color{black}\mE}$.
In this context, the number of travellers of type $\ell\in{\color{black}\mL}$ awaiting in the system at time $t$ is $ \sum_{e\in{\color{black}\mE}}N_{\ell,e}(t)$.
We consider the case where both drivers and riders will queue in the system.
An example about a bipartite-graph based on the matching mechanism is provided in \ref{app:example:2} of the supplementary material.

When a traveler of type $\ell\in{\color{black}\mL}$ arrives at time $t\geq 0$ and the system controller decides to serve it with a match $e\in\mE_{\ell}$, which corresponds to the edge $(\ell,\ell')$ or $(\ell',\ell)\in\mE$,
\begin{enumerate}[label=(\roman*)]
\item if there is a traveler of type $\ell'$ waiting there (that is, $N_{\ell',e}(t)>0$),  then both travelers will leave with a shared vehicle, and $N_{\ell',e}(t)$ decrements by one;  \label{case:1}
\item if there is no traveler of type $\ell'$ waiting there (that is, $N_{\ell',e}(t)=0$), then the number of travelers of type $\ell$ awaiting to be served by the match $e$, $N_{\ell,e}(t)$, increments by one. \label{case:2}
\end{enumerate}
Note that in the \ref{case:1} case, although the newly arrived type-$\ell$ traveller will incur an increment in $N_{\ell,e}(t)$, this new type-$\ell$ traveller will leave with a type-$\ell'$ traveller immediately and simultaneously, causing a decrement in both $N_{\ell,e}(t)$ and $N_{\ell',e}(t)$. 
It follows with unchanged $N_{\ell,e}(t)$ and $N_{\ell',e}(t)$ decremented by one.
If the shared ride is unsuccessful,  then, as described at the end of Section~\ref{subsec:ride-sharing}, the travelers may join the ride-sharing system again with different origins, destinations, travel preferences et cetera.

Apart from the arrivals, reneging events also incur decrements in the number of the waiting travelers.
Recall that, as defined in Section~\ref{subsec:ride-sharing}, the waiting time, for which a traveler of type $\ell\in{\color{black}\mL}$ becomes annoyed when awaiting to be served by a match $e\in\mE_{\ell}$, is exponentially distributed with
a rate dependent on the match $e$.
We refer to this rate as the reneging rate, represented by $\mu_{\ell,e}\in\mathbb{R}_0$. Once an $\ell$-traveler reneges while waiting to be served by the match $e$, the value of $N_{\ell,e}(t)$ decrements by one.

\vspace{-0.2cm}
\subsubsection{\color{black}State Variables}

In this context, for any $t\geq 0$ and $e\in{\color{black}\mE}$ incident to vertices $\ell,\ell'\in{\color{black}\mL}$, if $N_{\ell,e}(t)>0$ then $N_{\ell',e}(t)=0$. For travelers of the same type waiting for the same match $e$, they will be matched in a first-come-first-serve manner.
Define $N< +\infty$ as the maximal number of the waiting travelers. For any $\ell\in{\color{black}\mL}$ and $e\in{\color{black}\mE}$, we do not allow more than $N$ travelers of type $\ell$ waiting for the match $e$. It follows that $N_{\ell,e}(t)$ is no greater than $N$. This upper bound imposed on the traveler numbers guarantees the stability of the system and, based on Little's formula, issues an upper bound on the average waiting time of all the travelers. 
{\color{black}To simplify notation, for $e\in\mE$ and its incident vertices $\ell \in\mL^D$ and $\ell'\in \mL\backslash \mL^D$, define \vspace{-0.5cm}
\begin{equation}
    N_e(t) = \begin{cases}
        N_{\ell,e}(t), & \text{if } N_{\ell,e}(t) \geq 0,\\
        -N_{\ell',e}(t), & \text{otherwise},
    \end{cases}
\end{equation} 
and omitting subscript $\ell$ throughout the remainder of this paper.
Meanwhile, for any $e\in\mE$ and $\ell\in\mL$, given $N_e(t)$, we uniquely obtain $N_{\ell,e}(t)$: if $\ell = \ell^D(e)$, $N_{\ell,e}(t) = \max\{N_e(t),0\}$; if $\ell=\ell^R(e)$, $N_{\ell,e}(t) = \max\{-N_e(t),0\}$; otherwise, $N_{\ell,e}(t)=0$. }
Define $\tilde{\mathscr{N}}\coloneqq {\color{black}\{-N,-N+1,\ldots,N\}}$ as the set of possible values of {\color{black}$N_e(t)$} for any $e\in{\color{black}\mE}$ and $t\geq 0$.
If the total number of waiting travelers of type $\ell\in{\color{black}\mL^D}$ {\color{black}(or $\mL\backslash \mL^D$)} at time $t\geq 0$, $\sum_{e\in\mE_{\ell}}N_e(t)$ {\color{black}(or $-\sum_{e\in\mE_{\ell}}N_e(t)$)}, reaches the upper bound $N |\mE_{\ell}|$, then a newly arrived traveler of type $\ell$ will be rejected. The probability of the rejection decreases to zero as $N$ tends to infinity. 
To address the rejection case, for each traveler type $\ell\in{\color{black}\mL}$, define $e(\ell)\coloneqq E+\ell$ as a \emph{virtual match} that is associated with a stochastic process $\{N_{e(\ell)}(t), t\geq 0\}$ where $N_{e_{\ell}}(t) \equiv 0$ for all $t\geq 0$ and this process incurs no reward or cost. 
If a traveler of type $\ell$ is rejected, then we alternatively say that this traveler is assigned to the virtual match $e(\ell)$.
Let {${\color{black}n}(t)\coloneqq \Bigl({\color{black}N_e(t)}:\ e\in{\color{black}\mE}\Bigr)$} represent the system state 
at time $t$, and let
$\mathscr{N}\coloneqq \tilde{\mathscr{N}}^E$
represent the state space of the counting process $\{{\color{black}n}(t),t\geq0\}$.
\vspace{-0.2cm}
\subsubsection{\color{black} Action Variables}

The system controller makes decisions upon arrivals of travellers.
Consider state-dependent action variables $a_{\ell,e}({\color{black}n}(t))\in\{0,1\}$, for $t\geq 0$, $\ell\in{\color{black}\mL}$ and $e\in\mE_{\ell}\cup \{e(\ell)\}$.
For $\ell\in{\color{black}\mL}$, $e\in\mE_{\ell}\cup\{e(\ell)\}$ and a new type-$\ell$ traveler arrives at $t\geq 0$, if $a_{\ell,e}({\color{black}n}(t))=1$, then the newly arrived traveler of type $\ell$ is assigned to the match $e$; otherwise, the traveler will not be served by the match $e$. 
In the former case, the newly arrived type-$\ell$ traveler will be served by the match $e$. If this match $e$ is an eligible match $e\in\mE_{\ell}$, this assignment incurs state transitions in line with the manner described in \ref{case:1} and \ref{case:2}, right after its arrival in time $t$.
In this context, these action variables should satisfy \vspace{-0.5cm}
\begin{equation}\label{eqn:action_constraint}
\sum\nolimits_{e\in\mE_{\ell}\cup\{e(\ell)\}} a_{\ell,e}({\color{black}n}(t)) = 1,~\forall t \geq 0,~\ell\in{\color{black}\mL}, \vspace{-0.3cm}
\end{equation}
which guarantees that, at time $t\geq 0$, a newly arrived traveler of type $\ell\in{\color{black}\mL}$ will be assigned to one eligible match $e\in\mE_{\ell}$ or be rejected by selecting $e(\ell)$. 
Also, as we do not allow more than $N$ travelers waiting there, let \vspace{-0.5cm}
\begin{equation}\label{eqn:uncontrollable}
a_{\ell,e}({\color{black}n}(t))\mathds{1}\{N_{\ell,e}(t)=N\} = 0,~\forall t\geq 0, \ell\in{\color{black}\mL},e\in\mE_{\ell}, \vspace{-0.3cm}
\end{equation}
where $\mathds{1}$ is the indicator function, {\color{black}and $N_{\ell,e}(t)$ is uniquely determined with given $N_e(t)$}.
Let $\bm{a}({\color{black}n}(t)) \coloneqq (a_{\ell,e}({\color{black}n}(t)):\ell\in{\color{black}\mL},e\in \mE_{\ell})$.
We do not include the $a_{\ell,e(\ell)}({\color{black}n}(t))$ in the action vector $\bm{a}({\color{black}n}(t))$ because, based on \eqref{eqn:action_constraint}, it can {\color{black} be induced} by the values of $a_{\ell,e}({\color{black}n}(t))$ for all $e\in\mE_{\ell}$.
In this context, for $e\in{\color{black}\mE}$, the state transition rates of the underlying process $\bigl\{{\color{black}N_e(t)}, t\geq 0\bigr\}$ are presented in Figure~\ref{fig:markov_chain} in \ref{app:example:2} of the supplementary material.

Here, the transition rates between states are affected by arrivals, departures with shared rides, and reneges of travelers.
As described in \ref{case:1} and \ref{case:2}, upon a traveler arrival, based on an employed scheduling policy, the traveler will be assigned to a specific eligible match, incurring increment or decrement (caused by a departure with shared ride) of the associated state variable. 
As shown in Figure~\ref{fig:markov_chain} in \ref{app:example:2} of the supplementary material, the state variable increments are only incurred by arrivals of the corresponding travelers.
While, for the state variable decrements, apart from the traveler departures with shared rides (arrivals of the matched partners), 
the traveler reneging events also make contributions.
Recall that the waiting time, for which a traveler of type $\ell\in{\color{black}\mL}$ becomes annoyed when awaiting to be served by a match $e\in\mE_{\ell}$, is exponentially distributed with
the rate $\mu_{\ell,e}$.
When there are $N_{\ell,e}(t)=n$ travelers of type $\ell$ who are waiting for the match $e$, the remaining time until one of the $n$ travelers reneges is exponentially distributed with the rate $n\mu_{\ell,e}$.
Together with the departures of shared rides (arrivals of the matched partners), when {\color{black}$N_e(t)=n\geq 1$},  the rate for $N_{\ell^D(e),e}(t)$ transitioning from $n$ to $n-1$ is $a_{\ell^R(e),e}({\color{black}n}(t))\lambda_{\ell^R(e)} + n\mu_{\ell^D(e),e}$, where the action variable $a_{\ell^R(e),e}({\color{black}n}(t))$ decides whether the incoming $\ell^R(e)$-riders will be assigned to match the $\ell^D(e)$-drivers.  
{\color{black}For $N_e(t) = -n\leq -1$ transitioning to $-n+1$, the process $\{N_e(t),t\geq 0\}$ evolves in a similar manner.}

The action variables $\bm{a}({\color{black}n}(t))$ for all $t\geq 0$ determine a scheduling policy $\phi$ for the process $\{{\color{black}n}(t),t\geq 0\}$.
Let $\Phi$ denote the set of all such policies determined by $\bm{a}({\color{black}n}(t))\in\{0,1\}^{\sum_{\ell\in{\color{black}\mL}}|\mE_{\ell}|}$. 
To highlight the relationship between the stochastic process and its underlying policy, we add superscript $\phi\in\Phi$ to its state and action variables. We rewrite ${\color{black}n}(t)$ and $\bm{a}({\color{black}n}(t))$ as ${\color{black}n}^{\phi}(t)$ and $\bm{a}^{\phi}({\color{black}n}^{\phi}(t))$, respectively.
The policies $\phi\in\Phi$ discussed in this paper are non-anticipative.
Although, mathematically, arrivals, departures and {\color{black}reneging} of travelers will all incur state transitions and the action vector $\bm{a}^{\phi}({\color{black}n}^{\phi}(t))$ potentially changes its value whenever the system state ${\color{black}n}^{\phi}(t)$ becomes different, 
we only need to read the value of the action vector upon new arrivals of travellers.

\vspace{-0.2cm}
\subsubsection{\color{black} Stochastic Optimisation}

At any time $t\geq 0$, if a rider and a driver depart together with a shared vehicle with respect to the match $e\in{\color{black}\mE}$, the ride-matching system earns a random reward with expectation $R_e\in\mathbb{R}_+$.
Whenever a traveler of type $\ell\in{\color{black}\mL}$ waiting to take a ride through the match $e\in{\color{black}\mE}$ becomes impatient and leaves alone, the ride-matching system is penalized by $C_{\ell,e}\in\mathbb{R}_0$. 
Alternatively, we say that the system generates $-C_{\ell,e}$ reward.
We aim to maximize the long-run average total reward of the process $\{{\color{black}n}^{\phi}(t), t\geq 0\}$, where the total reward is the sum of the rewards generated by serving all the travelers.
For $e\in{\color{black}\mE}$ and $t\geq 0$, given {\color{black}$N_e^\phi(t) = n\in \tilde{\mathscr{N}}$}, $a^{\phi}_{\ell^D(e),e}({\color{black}n}^{\phi}(t))=a_1$ and $a^{\phi}_{\ell^R(e),e}({\color{black}n}^{\phi}(t))=a_2$, the expected reward per time unit generated by the process $\bigl\{{\color{black}N_e^\phi(t)},t\geq 0\bigr\}$ is \vspace{-0.5cm}
\begin{subequations}
	\label{eqn:reward_func}
	\begin{empheq}[left={R_e({\color{black}n},a_1,a_2)\coloneqq\empheqlbrace\,}]{align}
	& R_e a_1\lambda_{\ell^D(e)} -C_{\ell^R(e),e} {\color{black}|n|} \mu_{\ell^R(e),e}, \text{ if } {\color{black}n<0},
	\label{eqn:reward_func_1} \\
	& R_e a_2\lambda_{\ell^R(e)} -C_{\ell^D(e),e} {\color{black}|n|} \mu_{\ell^D(e),e}, \text{ if } {\color{black}n>0},
	\label{eqn:reward_func_2}\\
	& 0,  \hspace{5.1cm}\text{otherwise}.
	\label{eqn:reward_func_3}
	\end{empheq} \vspace{-0.5cm}
\end{subequations}

We refer to $R_e({\color{black}n},a_1,a_2)$ as the \emph{reward function} of the process $\bigl\{{\color{black}N_e^\phi(t)}, t\geq 0\bigr\}$.
The reward function takes into account the three possible scenarios as demonstrated in \eqref{eqn:reward_func_1}-\eqref{eqn:reward_func_3}. 
In the first case with ${\color{black}n<0}$ in \eqref{eqn:reward_func_1}, 
$\lambda_{\ell^D(e)}$ is the arrival rate of the drivers of type $\ell^D(e)$ and, if $a_1=1$, it is also the rate for a rider of type $\ell^R(e)$ to be matched with a newly arrived driver of type $\ell^D(e)$.
In \eqref{eqn:reward_func_1}, $R_e a_1\lambda_{\ell^D(e)}$ represents the expected reward rate at time $t$ that is generated by matching arrived drivers of type $\ell^D(e)$ to the riders of type $\ell^R(e)$. 
For the second term in \eqref{eqn:reward_func_1}, since ${\color{black}|n|}\mu_{\ell^R(e),e}$ is the reneging rate of the waiting riders of type $\ell^R(e)$, $-C_{\ell^R(e),e}{\color{black}|n|}\mu_{\ell^R(e),e}$ is the expected reward rate at time $t$ generated by impatient riders waiting for the match $e$.
Equation~\eqref{eqn:reward_func_2} is formulated in the same vein.
In  \eqref{eqn:reward_func_3}, the system receives no reward.
We aim to maximize the total average reward
\vspace{-0.3cm}
\begin{equation}\label{eqn:objective}
\max\limits_{\phi\in\Phi} \lim\limits_{T\rightarrow +\infty}\frac{1}{T} \mathbb{E}\int_0^T\sum\limits_{e\in{\color{black}\mE}} R_e\Bigl({\color{black}N_e^\phi(t)},a^{\phi}_{\ell^D(e),e}({\color{black}n}^{\phi}(t)),a^{\phi}_{\ell^R(e),e}({\color{black}n}^{\phi}(t))\Bigr) d t, \vspace{-0.3cm}
\end{equation}
subject to \eqref{eqn:action_constraint} and \eqref{eqn:uncontrollable}.
We refer to this problem described by \eqref{eqn:objective}, \eqref{eqn:action_constraint} and \eqref{eqn:uncontrollable} as the \emph{ride-matching problem} (RMP). 
The process $\{{\color{black}n}^{\phi}(t), t\geq 0\}$ with the action variable $\bm{a}^{\phi}({\color{black}n}^{\phi}(t))$ is a continuous-time MDP that aims to maximize the long-run average reward described in \eqref{eqn:objective} with the state-dependent action variables subject to \eqref{eqn:action_constraint} and \eqref{eqn:uncontrollable}.

The RMP is complicated by the large state-space size which increases exponentially in the number of the matches $E=|\mE|$. 
For large $E$, optimal solutions become computationally intractable with conventional techniques, such as value iteration, and we resort to simple, heuristic scheduling policies that approximate the optimality.

The underlying stochastic process of our problem is a high-dimensional process consisting of parallel \emph{sub-processes}, $\bigl\{{\color{black}N_e^\phi(t)},t\geq 0\bigr\}$, $e\in {\color{black}\mE}$, which are coupled by the constraints of action variables stated in \eqref{eqn:action_constraint}.  
At each decision-making epoch, these constraints force at most one sub-process $e\in|\mE|$ to be selected for each traveler type $\ell\in{\color{black}\mL}$. A newly arrived traveler will be assigned to the selected match $e\in {\color{black}\mE}$.
Such a problem falls in the same vein of the RMABP, where only part of parallel stochastic sub-processes are/is selected.

The Whittle index policy proposed by \cite{whittle1988restless} selects sub-process(es) according to the descending order of state-dependent \emph{indices} mapped to these sub-processes.
These indices, referred to as the \emph{Whittle indices}, are
obtained by decomposing the original RMAB problem through the Whittle relaxation technique. 
More precisely, \cite{whittle1988restless} conjectured that if the optimal solution of a relaxed version of the RMAB problem exists in a threshold form, then the Whittle index policy, with the indices equal to the optimal thresholds, approaches optimality as the number of bandit/sub-processes tends to infinity. The existence of a threshold-form optimal solution to the relaxed version of the RMAB problem is referred to as the \emph{Whittle indexability}.
The indices are
computed independently through each sub-process, and the Whittle index policy exhibits significantly reduced computational complexity.


However, unlike the conventional RMAB problem (such as \cite{whittle1988restless,nino2007dynamic,graczova2014generalized,fu2021large}) or the extended RMAB problems with multi-actions for each bandit/sub-process (such as \cite{glazebrook2011general,hodge2015asymptotic,zayascaban2019asymptotically,killian2021qlearning,xiong2022reinforcement}), for $e\in{\color{black}\mE}$ and a given state ${\color{black}n}^{\phi}(t)={\color{black}n}$, the action variables $(a^{\phi}_{\ell^D(e),e}({\color{black}n}),a^{\phi}_{\ell^R(e),e}({\color{black}n}))$ are constrained by two different constraints in \eqref{eqn:action_constraint}: one for $\ell=\ell^D(e)$ and the other for $\ell=\ell^R(e)$.
It imposes two different Lagrange multipliers for each sub-process if the original RMP is relaxed in line with the Whittle relaxation technique.
Because of the two Lagrange multipliers that mutually depend on each other, the optimal solution for the relaxed version of the RMP will not appear in the threshold form requested by  Whittle indexability.
Nor does there exist a Whittle index or any critical real number, related to each state of each bandit/sub-process, that directly represents the marginal cost/reward of taking a certain action. 
In \citep{fu2018restless}, multiple Lagrange multipliers that are bonded with a single bandit process were discussed in the restless-bandit-based (RBB) resource allocation problem. 
\cite{fu2018restless} studied only binary actions for each bandit process;
while, here, the RMP focuses on matching behaviors of two different travelers and requests four actions for the sub-process associated with $e\in{\color{black}\mE}$: assign/not assign the newly arrived driver/rider to the eligible match $e$.
The non-binary actions for each sub-process further prevent conventional analysis on Whittle indexability or algorithms for numerically computing the indices from being applied here.


 \vspace{-0.5cm}
\section{Bivariate Index Policy}\label{sec:policies}
\vspace{-0.2cm}

The RMAB problem, in general, is proved in  \citep{papadimitriou1999complexity} to be PSPACE-complete, indicating that it is also NP-complete (see \citep{sipser1996introduction}). 
The challenges primarily arise from the inter-dependencies among the restless bandit/sub-processes that are coupled by the action constraint, which suggests a single Lagrange multiplier in the relaxed version of the RMAB problem.
It hints at unfavorable conditions for solving the RMP, which imposes stronger dependencies among various sub-processes.

\vspace{-0.4cm}
\subsection{Bivariate Indices} \label{subsec:duplex_index}


By randomizing action variables, we relax the original constraint~\eqref{eqn:action_constraint} to \vspace{-0.5cm}
\begin{equation}\label{eqn:action_constraint:relax}
\lim\limits_{T\rightarrow +\infty} \frac{1}{T}\mathbb{E}\Bigl[\int_0^T \sum\nolimits_{e\in\mE_{\ell}\cup\{e(\ell)\}} a^{\phi}_{\ell,e}({\color{black}n}^{\phi}(t))\Bigr] d t =1,~\forall\ell\in{\color{black}\mL}.\vspace{-0.5cm}
\end{equation}
The constraint \eqref{eqn:action_constraint:relax} is constructed by taking expectations at both sides of the original constraint~\eqref{eqn:action_constraint}.
The optimization problem with the objective \eqref{eqn:objective} subject to \eqref{eqn:action_constraint:relax} is a relaxed version of the original RMP. 
With Lagrange multipliers $\bm{\eta}\in\mathbb{R}_0^L$, the Lagrangian function of the relaxed problem is\vspace{-0.5cm}
\begin{multline}\label{eqn:lagrangian_function:special_case}
g(\bm{\eta})\coloneqq \max\limits_{\phi\in\Phi}\lim\limits_{T\rightarrow +\infty} \frac{1}{T}\mathbb{E}\int_0^T
\sum\limits_{e\in{\color{black}\mE}} \biggl(R_e\Bigl({\color{black}N_e^\phi(t)},a^{\phi}_{\ell^D(e),e}({\color{black}n}^{\phi}(t)),a^{\phi}_{\ell^R(e),e}({\color{black}n}^{\phi}(t))\Bigr)
\\ -\eta_{\ell^D(e)} a^{\phi}_{\ell^D(e),e}({\color{black}n}^{\phi}(t))-\eta_{\ell^R(e)} a^{\phi}_{\ell^R(e),e}({\color{black}n}^{\phi}(t)) \biggr) d t \\
+\lim\limits_{T\rightarrow +\infty} \frac{1}{T}\mathbb{E}\int_0^T \sum\limits_{\ell\in{\color{black}\mL}}\Bigl(-\eta_{\ell}a^{\phi}_{\ell,e(\ell)}\bigl({\color{black}n}^{\phi}(t)\bigr)\Bigr)dt
+ \sum\nolimits_{\ell\in{\color{black}\mL}}\eta_{\ell},\vspace{-0.5cm}
\end{multline}
where $\bm{\eta} = (\eta_{\ell}:\ell\in{\color{black}\mL})$ are the Lagrange multipliers for the $L$ constraints in~\eqref{eqn:action_constraint:relax}.
By the Strong Law of Large Numbers for Continuous Time Markov Chains (see \cite[Theorem 45 in Chapter 4]{serfozo2009basics}), the Lagrange function $g(\bm{\eta})$, of which the first part is the long-run average reward earned at jumping times of the process $\{{\color{black}n}^{\phi}(t), t>0\}$ under given policy $\phi$, is\vspace{-0.5cm}
\begin{multline}\label{eqn:lagrangian_function}
g(\bm{\eta}) = \max\limits_{\phi\in\Phi}\sum\limits_{e\in {\color{black}\mE}} \sum\limits_{{\color{black}n}\in\mathscr{N}}\pi^{\phi}({\color{black}n}) \biggl(R_e\Bigl({\color{black}n_e},a^{\phi}_{\ell^D(e),e}({\color{black}n}),a^{\phi}_{\ell^R(e),e}({\color{black}n})\Bigr)\\-\eta_{\ell^D(e)} a^{\phi}_{\ell^D(e),e}({\color{black}n})-\eta_{\ell^R(e)} a^{\phi}_{\ell^R(e),e}({\color{black}n})\biggr)
+\sum\limits_{\ell\in {\color{black}\mL}} \sum\limits_{{\color{black}n}\in\mathscr{N}}\pi^{\phi}({\color{black}n})(-\eta_{\ell})a^{\phi}_{\ell,e(\ell)}({\color{black}n}))
+ \sum\nolimits_{\ell\in{\color{black}\mL}}\eta_{\ell}, \vspace{-0.5cm}
\end{multline}
where $\pi^{\phi}({\color{black}n})$ is the stationary probability that the process $\{{\color{black}n}^{\phi}(t), t \geq 0\}$ is in state ${\color{black}n}$ when the policy is $\phi$.

For $e\in{\color{black}\mE}$, {\color{black}$n\in\tilde{\mathscr{N}}$}, and $\phi\in\Phi$, define action variable $\alpha^{\phi}_{\ell,e}{\color{black}(n)}\coloneqq \mathbb{E}\Bigl[a^{\phi}_{\ell,e}({\color{black}\bm{N}}^{\phi}(t)) \Bigr|  {\color{black}N^{\phi}_e(t)=n}\Bigr]$, and for $\ell\in{\color{black}\mL}$, define action variable for the virtual match $e(\ell)$, $\alpha^{\phi}_{\ell,e(\ell)} \coloneqq \mathbb{E}\Bigl[a^{\phi}_{\ell,e(\ell)}({\color{black}\bm{N}}^{\phi}(t))\Bigr]$,
taking values in $[0,1]$.
For $e\in{\color{black}\mE}$ and ${\color{black}n}\in\tilde{\mathscr{N}}$, let $\bm{\alpha}^{\phi}_e({\color{black}n})\coloneqq\Bigl( \alpha^{\phi}_{\ell,e}({\color{black}n}): \ell\in\{\ell^D(e),\ell^R(e)\}\Bigr)$.
Action variables $\bm{\alpha}^{\phi}_e({\color{black}n})$ for all ${\color{black}n}\in\tilde{\mathscr{N}}$ are sufficient to determine the expected reward rates and transition probabilities in all states of the
process $\bigl\{{\color{black}N_e^\phi(t)}, t\geq 0\bigr\}$.
{For $e\in{\color{black}\mE}$, define $\tilde{\Phi}_e$ as the set of all policies $\phi$ determined by $\bm{\alpha}^{\phi}_e({\color{black}n})$ for all ${\color{black}n}\in \tilde{\mathscr{N}}$, and}
in the maximization problem on the right hand side of \eqref{eqn:lagrangian_function}, there is no constraint
that restricts the maximal reward rate of the sub-process $\bigl\{{\color{black}N_e^\phi(t)}, t\geq 0\bigr\}$ for any $e\in{\color{black}\mE}$ once the others are known. 
In the same vein as the relaxation technique presented in \citep{whittle1988restless}, we obtain \vspace{-0.5cm}
\begin{multline}\label{eqn:sub_problems}
g(\bm{\eta}) = \sum\limits_{e\in {\color{black}\mE}} \max\limits_{\phi\in\tilde{\Phi}_e}\sum\limits_{{\color{black}n}\in\tilde{\mathscr{N}}}\pi^{\phi}_e({\color{black}n}) \biggl(R_e\Bigl({\color{black}n},\alpha^{\phi}_{\ell^D(e),e}({\color{black}n}),\alpha^{\phi}_{\ell^R(e),e}({\color{black}n})\Bigr)\\-\eta_{\ell^D(e)} \alpha^{\phi}_{\ell^D(e),e}({\color{black}n})-\eta_{\ell^R(e)} \alpha^{\phi}_{\ell^R(e),e}({\color{black}n})\biggr) 
+\sum\limits_{\ell\in {\color{black}\mL}} \max\limits_{\alpha^{\phi}_{\ell,e(\ell)}\in[0,1]}(-\eta_{\ell})\alpha^{\phi}_{\ell,e(\ell)}
+\sum\limits_{\ell\in{\color{black}\mL}}\eta_{\ell}, \vspace{-0.5cm}
\end{multline}
where $\pi^{\phi}_e({\color{black}n})$ is the stationary probability that the process $\bigl\{{\color{black} N_e^\phi}(t),t\geq 0\bigr\}$ stays in state ${\color{black}n}$.
The optimization on the right-hand side of \eqref{eqn:sub_problems} can be decomposed into $E+L$ independent sub-problems, corresponding to the $E$ matches in the first term and the $L$ virtual matches in the second term of \eqref{eqn:sub_problems}. 

{Consider the maximization for $e\in{\color{black}\mE}$ in the first item of \eqref{eqn:sub_problems}.}
For $e\in{\color{black}\mE}$, $\bm{\eta}\in\mathbb{R}^2$, and ${\color{black}n}\in\tilde{\mathscr{N}}\backslash\{{\color{black}0}\}$, consider the value function $V^{\bm{\eta}}_e({\color{black}n})$ that solves \vspace{-0.3cm}
\begin{equation}\label{eqn:bellman_duplex}
V^{\bm{\eta}}_e({\color{black}n}) = \max\limits_{\bm{a}\in \mathcal{A}({\color{black}n})}\Bigl\{
\frac{R_e({\color{black}n},a_1,a_2)-a_1\eta_1-a_2\eta_2-g }{u_e({\color{black}n},\bm{a})} 
+\sum\nolimits_{{\color{black}n'}\in \tilde{\mathscr{N}}}p_e^{\bm{a}}({\color{black}n},{\color{black}n}')V^{\bm{\eta}}_e({\color{black}n'})
\Bigr\},\vspace{-0.3cm}
\end{equation}
where 
$\mathcal{A}({\color{black}n}) \coloneqq 
\bigl\{\{0\}\cup\{\mathds{1}\{{\color{black}n}<N\}\}\bigr\}\times\bigl\{\{0\}\cup\{\mathds{1}\{{\color{black}n>-N}\}\}\bigr\}$,
$g\in \mathbb{R}$ is an attached criterion for maximizing the average reward of the underlying process, 
and $1/u_e({\color{black}n},\bm{a})$ and $p_e^{\bm{a}}({\color{black}n},{\color{black}n'})$ are the expected sojourn time in state ${\color{black}n}$ and the transition probability from state ${\color{black}n}$ to ${\color{black}n'}$, respectively, under policy $\phi$ with action variable $\bm{\alpha}^{\phi}_e({\color{black}n})=\bm{a}\in\mathcal{A}({\color{black}n})$.
Define $V^{\bm{\eta}}_e({\color{black}0})\equiv 0$ for all $\bm{\eta}\in\mathbb{R}^2_0$ and $e\in{\color{black}\mE}$.

\begin{figure}[t]
\centering
\includegraphics[width=0.33\linewidth]{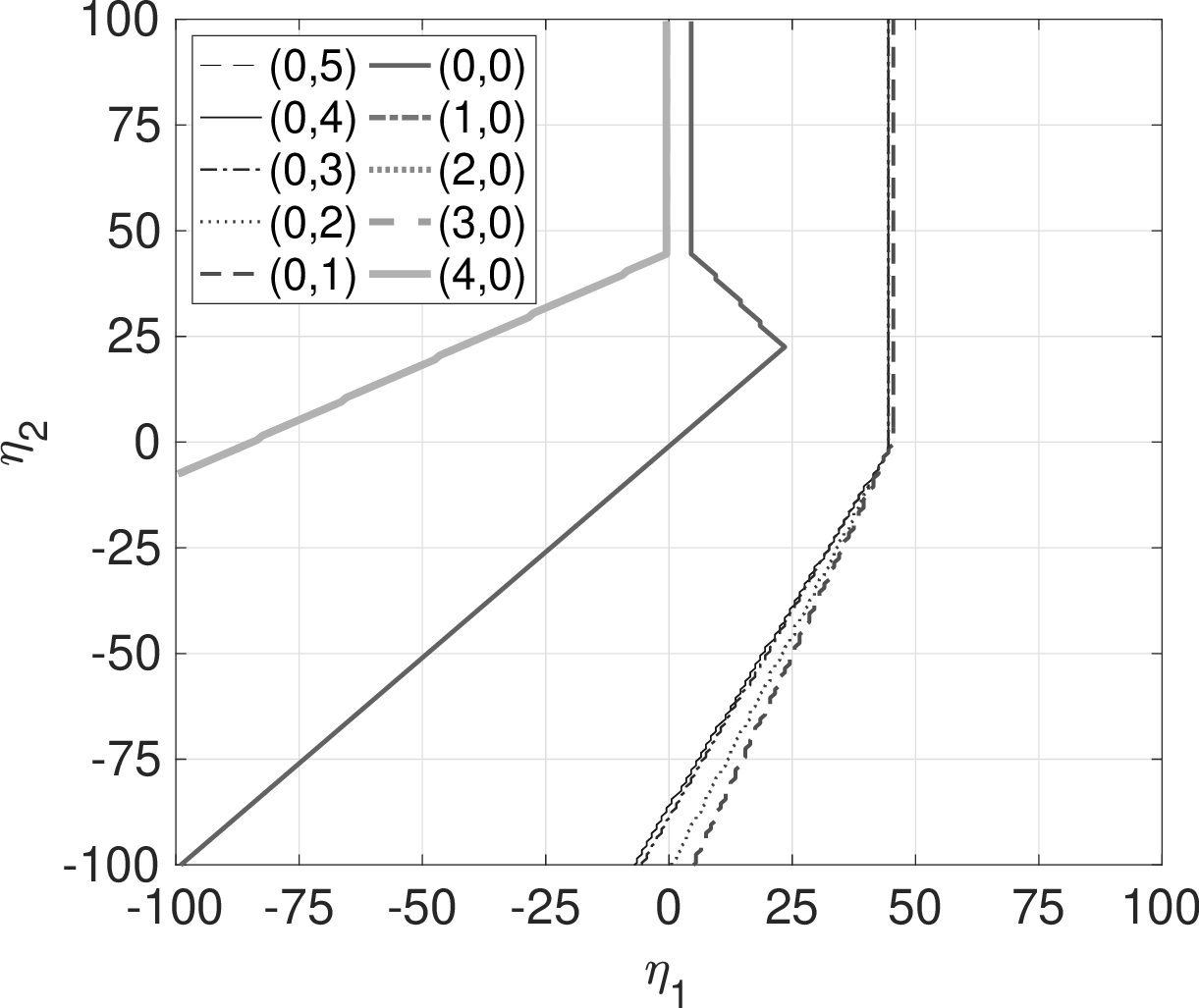}
\vspace{-0.5cm}\caption{Critical values of the multiplier vector $\bm{\eta}=(\eta_1,\eta_2)$ such that $a^{\phi}_{\ell^D(e),e}({\color{black}n})=0$ and $1$ are both optimal to the $(\bm{\eta},e)$-process, for different ${\color{black}n}\in\tilde{\mathscr{N}}\backslash\{{\color{black}N}\}$.
}\label{fig:example_1}
\vspace{-0.5cm}
\end{figure}

Equation~\eqref{eqn:bellman_duplex} is the Bellman equation for the underlying MDP of the process $\bigl\{{\color{black}N_e^\phi(t)}, t\geq 0\bigr\}$, of which the expected reward rate in state ${\color{black}n}$ under action $\bm{a}=(a_1,a_2)$ is $R_e({\color{black}n},a_1,a_2)-a_1\eta_1-a_2\eta_2$. We refer to such an MDP as an $(\bm{\eta},e)$-\emph{process} associated with the eligible match $e\in{\color{black}\mE}$ and aim to maximize its long-run average reward.
For ${\color{black}n}\in\tilde{\mathscr{N}}\backslash\{{\color{black}0}\}$, the value function $V^{\bm{\eta}}_e({\color{black}n})$ represents the maximal expected cumulative reward of a period of the $(\bm{\eta},e)$-process, which starts from state ${\color{black}n}$ and ends when it enters state ${\color{black}0}$.
In this context, for a given real number $g$ equal to the maximized average reward of the $(\bm{\eta},e)$-process, the actions $\bm{\alpha}^{\phi}_e({\color{black}n})=\bm{a}$ maximizing the right hand side of \eqref{eqn:bellman_duplex} for all states ${\color{black}n}\in\tilde{\mathscr{N}}\backslash\{{\color{black}0}\}$ form an optimal solution $\phi$ of the $(\bm{\eta},e)$-process.
The optimal action vector $\bm{\alpha}^{\phi}_e({\color{black}0})$ for state ${\color{black}0}$ with respect to the average reward of the $(\bm{\eta},e)$-process is obtained by maximizing\vspace{-0.3cm}
\begin{equation}\label{eqn:bellman_zero}
H^{\bm{\eta}}_e= \max\biggl\{0,\max\limits_{\bm{a}\in\mathcal{A}({\color{black}0})\backslash\{{\color{black}0}\}}\Bigl\{
\frac{R_e({\color{black}0},a_1,a_2)-a_1\eta_1-a_2\eta_2-g }{u_e({\color{black}0},\bm{a})} 
+\!\!\sum_{{\color{black}n'}\in \tilde{\mathscr{N}}}\!p_e^{\bm{a}}({\color{black}0},{\color{black}n'})V^{\bm{\eta}}_e({\color{black}n'})
\Bigr\}\biggr\},\vspace{-0.3cm}
\end{equation}
where $g$ is given and equal to the maximized average reward of the $(\bm{\eta},e)$-process.
If $H_e^{\bm{\eta}}=0$, then the optimal action vector $\bm{\alpha}^{\phi}_e({\color{black}0})=\bm{0}$; otherwise, it is equal to the action vector $\bm{a}$ that maximizes the second argument in \eqref{eqn:bellman_zero}.
Note that, for given $\bm{\eta}$, the optimal solution of the $(\bm{\eta},e)$-process with respect to the average reward may not be unique. 
Let $\Phi^*_e(\bm{\eta})$ ($e\in{\color{black}\mE}$, $\bm{\eta}\in \mathbb{R}^2$) represent the set of such optimal policies associated with the match $e$.

Observing \eqref{eqn:bellman_duplex} and \eqref{eqn:bellman_zero}, for each $e\in{\color{black}\mE}$ and $\bm{\eta}\in\mathbb{R}^2$, both multipliers in $\bm{\eta}$ play important roles in the achieved optimum of the $(\bm{\eta},e)$-process.
Recall that, unlike the conventional (multi-action) RMAB problems, here, the optimum is dependent on two different multipliers $\bm{\eta}\in\mathbb{R}^2$.
All the past theorems are no longer applicable when the dimension of $\bm{\eta}$ is greater than one. 
As an example, we plot in Figure~\ref{fig:example_1} critical values of the vector $\bm{\eta}=(\eta_1,\eta_2)$ such that $a^{\phi}_{\ell^D(e),e}({\color{black}n})=0$ and $1$ are both optimal to the $(\bm{\eta},e)$-process. 
Detailed simulation settings and maximization algorithms are provided in \ref{app:settings:example_1} of the supplementary material. 
Recall that, as defined in Section~\ref{subsec:graph}, for the eligible match $e\in{\color{black}\mE}$, the state space for the $(\bm{\eta},e)$-process is $\tilde{\mathscr{N}} = \{{\color{black}-N, -N+1,\dots,N}\}$, and $\alpha^{\phi}_{\ell,e}(N)\equiv 0$ under the optimality assumption.

In Figure~\ref{fig:example_1}, each curve corresponds to a state ${\color{black}n}\in\tilde{\mathscr{N}}\backslash\{{\color{black}N}\}$.
It is a set of specific points $\bm{\eta}=(\eta_1,\eta_2)$,
for which $a_1=0$ and $1$ both maximize the right hand side of \eqref{eqn:bellman_duplex} for the state ${\color{black}n}$, when the given parameter $g$ is equal to the maximized average reward of the $(\bm{\eta},e)$-process. That is, there exist two optimal policies $\phi_1,\phi_2\in\Phi^*_e(\bm{\eta})$ with $\alpha^{\phi_1}_{\ell^D(e),e}({\color{black}n})=1$ and $\alpha^{\phi_2}_{\ell^D(e),e}({\color{black}n})=0$, respectively, such that both $\phi_1$ and $\phi_2$ maximize the average reward of the $(\bm{\eta},e)$-process.

In Figure~\ref{fig:example_1}, each curve divides the plane into two areas: left and right ones.
In our simulations, there exists an optimal policy $\phi$ to the $(\bm{\eta},e)$-process such that, if $\bm{\eta}$ locates at the left area of a curve for state ${\color{black}n}\in\tilde{\mathscr{N}}\backslash\{{\color{black}N}\}$, then $\alpha^{\phi}_{\ell^D(e),e}({\color{black}n})=1$; if it locates at the right area, $\alpha^{\phi}_{\ell^D(e),e}({\color{black}n})=0$; otherwise, both $\alpha^{\phi}_{\ell^D(e),e}({\color{black}n})= 1$ and $0$ are optimal.
For ${\color{black}n}\in\tilde{\mathscr{N}}\backslash\{{\color{black}N}\}$, $e\in{\color{black}\mE}$, and a policy $\phi$ optimal to the $(\bm{\eta},e)$-process, we refer to the sets of $\bm{\eta}\in\mathbb{R}^2$ with $\alpha^{\phi}_{\ell^D(e),e}({\color{black}n})=1$ and $0$ as \emph{active} and \emph{passive areas}, respectively, of the \emph{first action variable} in the \emph{state-match} pair (SM pair) $({\color{black}n},e)$.

Unlike the conventional RMAB with a one-dimensional Lagrange multiplier for each sub-problem, here, we are not able to identify a critical real number for each SM pair, which represents the marginal reward of this SM pair and serves as a comparable parameter (for instance, the Whittle index) that quantifies the priority of the SM pair.
In the RMP case with two-dimensional $\bm{\eta}$, it necessitates an interpretation of the two-dimensional multipliers and, more importantly, a discussion on how this interpretation can help construct an appropriate scheduling strategy applicable to the original RMP.

Let $\mathcal{G}$ represent the set of all real-valued functions of $\eta_2\in\mathbb{R}$ that are continuous in $\eta_2$, and consider the following definitions.
For $e\in{\color{black}\mE}$, ${\color{black}n}\in\tilde{\mathscr{N}}$,
define a set of functions of $\bm{\eta}\in\mathbb{R}^2$\vspace{-0.5cm}
\begin{multline} 
\mathcal{B}_1({\color{black}n},e) \coloneqq 
\Bigl\{\map \in \mathcal{G}~\Bigl|~\text{For all }\bm{\eta}\in \mathbb{R}^2 \text{ with } \eta_1=\Gamma(\eta_2), \text{ there exist } \phi_1,\phi_2\in \Phi_e^*(\bm{\eta})\\ \text{ such that }\alpha^{\phi_1}_{\ell^D(e),e}({\color{black}n}) = 1,\alpha^{\phi_2}_{\ell^D(e),e}({\color{black}n}) = 0\Bigr\}. 
\vspace{-0.5cm}
\end{multline}
In particular, 
if ${\color{black}n=N}$ then $\mathcal{B}_1({\color{black}n},e)=\emptyset$.
For $e\in{\color{black}\mE}$, ${\color{black}n}\in\tilde{\mathscr{N}}$, and $\Gamma\in\mathcal{B}_1({\color{black}n},e)$, 
define  \vspace{-0.5cm}
\begin{equation}
\mathscr{A}^{\map}_1({\color{black}n},e) \coloneqq 
\Bigl\{\bm{\eta} \in \mathbb{R}^2~\Bigl|~\text{there exists } \phi\in \Phi_e^*(\bm{\eta})\text{ such that }~\alpha^{\phi}_{\ell^D(e),e}({\color{black}n}) = 1\text{ and } \eta_1 \neq\Gamma(\eta_2)\Bigr\}.\vspace{-0.5cm}
\end{equation}
For any function $\Gamma\in\mathscr{B}_1({\color{black}n},e)$ and all the points $\bm{\eta}\in\mathbb{R}^2$ satisfying $\eta_1=\Gamma(\eta_2)$, there is no difference to take $\alpha^{\phi}_{\ell^D(e),e}({\color{black}n})=1$ or $0$ - they are both optimal to the $(\bm{\eta},e)$-process.
The set $\mathscr{A}^{\Gamma}_1({\color{black}n},e)$ is a subset of the active area of the SM pair $({\color{black}n},e)$, excluding those $\bm{\eta}$ with $\eta_1=\Gamma(\eta_2)$. 
We refer to $\mathscr{A}^{\Gamma}_1({\color{black}n},e)$ as the \emph{active area with respect to  $\Gamma$} for SM pair $({\color{black}n},e)$.
For each SM pair $({\color{black}n},e)$ with ${\color{black}n}\in\tilde{\mathscr{N}}\backslash\{{\color{black}N}\}$ and $e\in{\color{black}\mE}$, there always exists $\Gamma\in \mathscr{B}_1({\color{black}n},e)$.

In the example presented in Figure~\ref{fig:example_1}, each of the curve corresponds to a state ${\color{black}n}\in\tilde{\mathscr{N}}\backslash\{{\color{black}N}\}$ and can be described by a  $\Gamma_{{\color{black}n}}\in\mathscr{B}_1({\color{black}n},e)$.
We observe that the active areas with respect to the $\Gamma_{{\color{black}n}}$ (${\color{black}n}\in\tilde{\mathscr{N}}\backslash\{{\color{black}N}\}$) are \emph{nested} with each other.
That is, there is a permutation of states in $\tilde{\mathscr{N}}\backslash\{{\color{black}N}\}$, denoted by $\{{\color{black}n}_1, {\color{black}n}_2,\ldots, {\color{black}n}_{2N}\}$, such that $\mathcal{A}^{\Gamma_{{\color{black}n}_1}}_1({\color{black}n}_1,e)\subset \mathcal{A}^{\Gamma_{{\color{black}n}_2}}_1({\color{black}n}_2,e)\subset \ldots \subset \mathcal{A}^{\Gamma_{{\color{black}n}_{2N}}}_1({\color{black}n}_{2N},e)$.
In the same vein as the Whittle indices for the conventional RMAB problem, the nesting phenomenon indicates \emph{priorities} of these states: for given $\bm{\eta}\in\mathbb{R}^2$, $i=1,2,\ldots,2N-1$, and the $(\bm{\eta},e)$-process, if $\bm{\eta}\in\mathcal{A}_1^{\Gamma_{{\color{black}n}_i}}({\color{black}n}_i,e)$ (that is, $\alpha^{\phi}_{1,e}({\color{black}n_i}) = 1$ is optimal), then $\alpha^{\phi}_{1,e}({\color{black}n_j})=1$ is also optimal for $j=i+1,i+2,\ldots,2N$, because $\mathcal{A}_1^{\Gamma_{{\color{black}n}_i}}({\color{black}n}_i,e)\subset \mathcal{A}_1^{\Gamma_{{\color{black}n}_j}}({\color{black}n}_j,e)$.
The state ${\color{black}n}_j$ is prioritized against state ${\color{black}n}_i$ for taking optimal action variable $\alpha^{\phi}_{\ell^D(e),e}({\color{black}n_j})=1$.

We can quantify the SM-pair priorities by comparing the values of the $\eta_1=\Gamma(\eta_2)$, where $\Gamma\in\mathscr{B}_1({\color{black}n},e)$, for a given $\eta_2\in\mathbb{R}$, because larger $\Gamma(\eta_2)$ implies a larger active area.
Unfortunately, for ${\color{black}n}\in\tilde{\mathscr{N}}\backslash\{(N,0)\}$ and $e\in{\color{black}\mE}$, although there always exists $\Gamma\in\mathscr{B}_1({\color{black}n},e)$, it is unclear how to compute this $\Gamma$ or obtain the exact value of $\eta_1=\Gamma(\eta_2)$ for any $\eta_2\in\mathbb{R}$.

To approximate the priorities quantified by $\Gamma\in\mathscr{B}_1({\color{black}n},e)$, we construct a surrogate problem of the $(\bm{\eta},e)$-process. 
Consider the process $\bigl\{{\color{black}N_e^\phi(t)}, t\geq 0\bigr\}$, which takes  $R_e({\color{black}n},a_1,1)$ as the expected reward rate in state ${\color{black}n}\in\tilde{\mathscr{N}}\backslash\{-N\}$ when the action variable $\alpha^{\phi}_{\ell^D(e),e}({\color{black}n}) = a_1$, and, for state ${\color{black}-N}$, the expected reward rate of the process is $R_e({\color{black}-N},a_1,0)$ with $\alpha^{\phi}_{\ell^D(e),e}({\color{black}n}) = a_1$. 
It follows with an MDP similar to the $(\bm{\eta},e)$-process when $\bm{\eta}=\bm{0}$, except that $\alpha^{\phi}_{\ell^R(e),e}({\color{black}n})\equiv 1$ for all ${\color{black}n}\in\tilde{\mathscr{N}}\backslash\{{\color{black}-N}\}$.
We refer to such an MDP as the $e$-process. 
The $e$-process is a standard restless bandit process, for which the Whittle indices are computable. 
Let $\eta_e^1({\color{black}n})$ represent the Whittle index of the state ${\color{black}n}\in\tilde{\mathscr{N}}\backslash\{{\color{black}N}\}$ of the $e$-process. 
Consider the following proposition.

\begin{proposition}\label{prop:priority}
For $e\in{\color{black}\mE}$ and ${\color{black}n}\in\tilde{\mathscr{N}}\backslash \{{\color{black}N}\}$,
if $\mu_{\ell^R(e)}=\mu_{\ell^D(e)}=0$ and $\lambda_{\ell^R(e)} = \lambda_{\ell^D(e)}$, then
there exist $E\in\mathbb{R}$ and $\Gamma\in\mathcal{B}_1({\color{black}n},e)$ such that, for any $\eta_2 < E$, $\Gamma(\eta_2) = \eta_e^1({\color{black}n})-\eta_2$.\vspace{-0.3cm}
\end{proposition}
The proof of Proposition~\ref{prop:priority} is provided in \ref{app:priority} of the supplementary material. 
{\color{black}For an insight of the proof, 
when $\eta_2$ is sufficiently small, it is optimal to have $\alpha^{\phi}_{\ell^R(e),e}(n)=1$ for all $\tilde{\mathscr{N}}\backslash\{-N\}$.
It follows with a $(\bm{\eta},e)$-process the same as the $e$-process except that, when $N^{\phi}_e(t)\in \tilde{\mathscr{N}}\backslash\{N,-N\}$, the $(\bm{\eta},e)$-process' reward rate is $R_e\Bigl(N^{\phi}_e(t),\alpha^{\phi}_{\ell^D(e),e}\bigl(N^{\phi}_e(t)\bigr),1\Bigr) - \alpha^{\phi}_{\ell^D(e),e}\bigl(N^{\phi}_e(t)\bigr)\eta_1$ other than $R_e\Bigl(N^{\phi}_e(t),\alpha^{\phi}_{\ell^D(e),e}\bigl(N^{\phi}_e(t)\bigr),1\Bigr)$ for the $e$-process.
Such an extra term $-\alpha^{\phi}_{\ell^D(e),e}\bigl(N^{\phi}_e(t)\bigr)\eta_1$ in general incurs deviations between $\Gamma^{n,e}(\eta_2)$ and $\eta^1_e(n)$ for a $\Gamma^{n,e}\in\mathscr{B}_1(n,e)$.
The hypothesis of Proposition~\ref{prop:priority} implies a balanced situation where the extra term equally affects $\Gamma^{n,e}(\eta_2)$ for all $n \in \tilde{\mathscr{N}}\backslash\{N\}$ and $e\in\mE$.
In this case, for any $n_1,n_2\in\tilde{\mathscr{N}}\backslash\{N\}$ and $e_1,e_2\in\mL$, $\Gamma^{n_1,e_1}(\eta_2)-\Gamma^{n_2,e_2}(\eta_2)$ becomes independent of $\eta_2$; that is, the order of $\Gamma^{n_1,e_1}(\eta_2)$ and $\Gamma^{n_2,e_2}(\eta_2)$ is independent of $\eta_2$.
}

Proposition~\ref{prop:priority} suggests that, when the reneging rates are negligible and arrival rates of the drivers and riders are compatible to each other, the order of $\Gamma(\eta_2)$ is decided by $\eta_e^1({\color{black}n})$ and independent from the value of $\eta_2$ - the SM priorities can be quantified by $\eta_e^1({\color{black}n})$. 
Higher $\eta_e^1({\color{black}n})$ implies a larger active area with respect to the $\Gamma$ and higher priority of the SM pair $({\color{black}n},e)$.
Recall that the $e$-process is a surrogate problem constructed to quantify the SM priorities and is not any real process discussed in this paper. 
The $\eta_e^1({\color{black}n})$ for each SM pair $({\color{black}n},e)$ is the Whittle index of state ${\color{black}n}$ of the $e$-process, which can be computed through conventional techniques, such as the bisection method in \citep{fu2019towards,fu2021large} and the marginal productivity (MP) index in \citep{nino2001restless,nino2002dynamic,nino2007dynamic,nino2020verification}, or exists in a closed form for a range of special cases in \citep{fu2023restless}.

For the more general case with significant reneging rates and unbalanced arrival rates, 
it is unclear whether or not the order of $\Gamma(\eta_2)$ mentioned in Proposition~\ref{prop:priority} can be independent of $\eta_2$. 
In Section~\ref{sec:example}, we numerically discuss the performance of the quantification through $\eta_e^1({\color{black}n})$.

Apart from the quantification for drivers discussed above, due to the symmetry of each $(\bm{\eta},e)$-process, we can also quantify the priorities of the SM pairs for the associated riders through a surrogate problem similar to the $e$-process.
Precisely, consider the process $\bigl\{{\color{black}N_e^\phi(t)}, t\geq 0\bigr\}$, for which the expected reward rate in state $({\color{black}n})\in\tilde{\mathscr{N}}\backslash\{{\color{black}N}\}$ when the action variable $\alpha^{\phi}_{\ell^R(e),e}({\color{black}n}) = a_2$ is $R_e({\color{black}n},1,a_2)$ and, for the edge state ${\color{black}N}$, the expected reward rate is $R_e({\color{black}N},0,a_2)$.
This is an MDP similar to the $e$-process, except that, here, we fix $\alpha^{\phi}_{\ell^D(e),e}({\color{black}n})\equiv 1$ (other than $\alpha^{\phi}_{\ell^R(e),e}({\color{black}n})$) for all $\tilde{\mathscr{N}}\backslash\{{\color{black}N}\}$. 
Such an MDP associated with $e\in{\color{black}\mE}$ is also a standard bandit process with computable Whittle indices, which is denoted by $\eta_e^2({\color{black}n})$ for state  ${\color{black}n}\in\tilde{\mathscr{N}}\backslash\{{\color{black}-N}\}$.
Due to the symmetry of the $(\bm{\eta},e)$-process, when we focus on the active areas of the action variable of the riders (that is, $\alpha^{\phi}_{\ell^R(e),e}({\color{black}n})$), the priority of the SM pair $({\color{black}n},e)$ can also be quantified by $\eta_e^2({\color{black}n})$.

We refer to $\eta_e^1({\color{black}n})$ and $\eta_e^2({\color{black}n})$ as the \emph{bivariate index} for the SM pair $({\color{black}n},e)$ with respect to the driver and rider, respectively.
The bivariate index for an SM pair $({\color{black}n},e)$ summarizes the variate of the two Lagrange multipliers, $\eta_1$ and $\eta_2$,  which intuitively represent the marginal rewards of selecting the match $e$ for the newly arrived driver and rider, respectively, when ${\color{black}N_e^\phi(t)} = {\color{black}n}$.


\vspace{-0.3cm}

\subsection{Bivariate Index Policy}\label{subsec:BI_policy}

For the RMP described in \eqref{eqn:objective} and \eqref{eqn:action_constraint}, we propose a policy that always selects the eligible matches with the highest bivariate indices upon the arrival of a driver or a rider. 
We refer to such a policy as the \emph{bivariate index policy} (BI policy).
More precisely, consider the stochastic process $\{{\bm{N}}^{\text{BI}}(t),t\geq 0\}$ under the BI policy with indices quantified by $\eta_e^1({\color{black}n})$ or $\eta_e^2({\color{black}n})$.
For $\ell\in{\color{black}\mL}$ and $e\in\mE_{\ell}$, the action variables of the BI policy are \vspace{-0.5cm}\begin{equation}\label{eqn:BI_policy}
a^{\text{BI}}_{\ell,e}(\bm{N}^{\text{BI}}(t))=
\begin{cases}
1,& \text{if } N^{\text{BI}}_{{\color{black}e}}(t) < N,~\ell\in {\color{black}\mL^D}, \text{ and }e = \arg\max\limits_{e'\in\mE_{\ell}}\eta_{e'}^1({\color{black}N^{\text{BI}}_{e'}}(t)),\\
1,& \text{if }N^{\text{BI}}_{\color{black}e}(t) > -N,~\ell\notin {\color{black}\mL^D}, \text{ and }e = \arg\max\limits_{e'\in\mE_{\ell}}\eta_{e'}^2({\color{black}{N}^{\text{BI}}_{e'}}(t)),\\
0, &\text{otherwise}, \vspace{-0.5cm}
\end{cases}
\end{equation}
where $a^{\text{BI}}_{\ell,e(\ell)}\bigl(\bm{N}^{\text{BI}}(t)\bigr) = 1 - \sum\nolimits_{e\in\mE_{\ell}} a^{\text{BI}}_{\ell,e}\bigl(\bm{N}^{\text{BI}}(t)\bigr)$.
If  $\arg\max$ returns more than one number, we break the tie by choosing the smallest.
We only care about the priorities of the SM pairs quantified by the bivariate indices.

The BI policy is applicable to the original RMP described in \eqref{eqn:objective} and \eqref{eqn:action_constraint}, for which the indices $\eta_e^1({\color{black}n})$ and $\eta_e^2({\color{black}n})$ are calculated offline with computational complexity linear to the total number of eligible matches $E$ and at most quadratic to the size of the state space of each sub-process $|\tilde{\mathscr{N}}|=2N+1$. 
Since the indices are computed independently per each $e\in{\color{black}\mE}$, they can also be obtained through parallel computing, which can further accelerate the computation process.
Because of the independence for computing indices,
the BI policy is directly applicable to a flexible system allowing newly added traveler types and/or eligible matches without demanding re-calculation of the indices for the old ones.
The storage complexity for the pre-computed indices is linear to $E$.
As defined in \eqref{eqn:BI_policy}, the online implementation of the BI policy is about sorting the matches according to their instantaneous indices $\eta_e^1\bigl({\color{black}N^{\text{BI}}_e(t)}\bigr)$ and $\eta_e^2\bigl({\color{black}N^{\text{BI}}_e(t)}\bigr)$, which can be instantiated by using max heaps. 
In this way, the computational complexity is logarithmic to $|\mE_{\ell}|$, where $\ell$ is the type of the newly arrived traveler.
The BI policy is scalable for a system with large $L$ and $E$.
{\color{black}We provide the pseudo-codes for implementing BI in \ref{app:whittle} of the supplementary material.}

\vspace{-0.6cm}
\section{Numerical Results}\label{sec:example}

\vspace{-0.3cm}
To demonstrate the effectiveness of the BI policy, we provide simulation results here, which compare the BI policy with baseline policies in a system based on realistic data estimated with time-varying ride demands from \citep{Chicago} and the local population from \citep{ABS2018}. 
Simulations for another system with time-invariant arrival rates and stationary performance are presented in \ref{app:simulation_uninetwork} of the supplementary material.


\begin{figure}[t]
\centering
\subfigure[]{\includegraphics[width=0.32\linewidth]{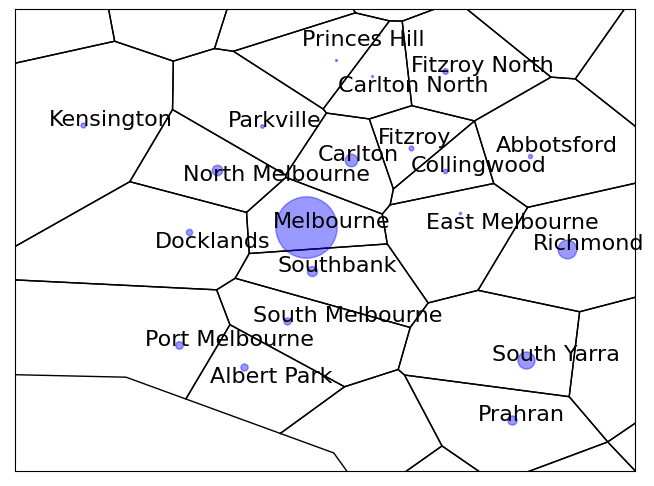}\label{fig:melbourne_example:demands}} 
\subfigure[]{\includegraphics[width=0.32\linewidth]{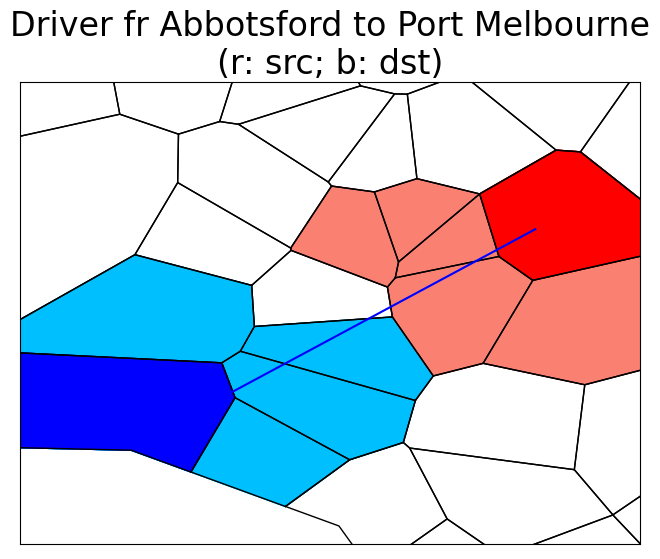}\label{fig:melbourne_example:2}}
\subfigure[]{\includegraphics[width=0.32\linewidth]{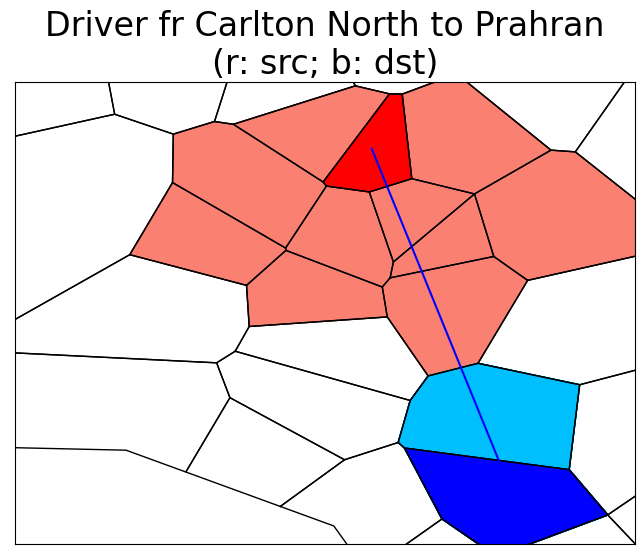}\label{fig:melbourne_example:3}}
\vspace{-0.5cm}\caption{Transportation network for Melbourne city: (a) districts within the 5km radius from Melbourne CBD, where traveler demands are represented by blue circles (the larger the circle, the higher the demand); (b) eligible traveler types that can match the traveler type $12$ (drivers from \emph{Abbotsford} to \emph{Port Melbourne}); and (c) eligible traveler types that can match the traveler type $70$ (drivers from \emph{Carlton North} to \emph{Prahran}).\vspace{-0.2cm}}\label{fig:melbourne_example}
\end{figure}

Consider the baseline policies: Join-in-the-Longest-Queue (JLQ), and a policy that myopically assigns travelers to the matches with the highest immediate expected return. 
For the match $e\in{\color{black}\mE}$ associated with traveler types $\ell^D(e),\ell^R(e)$, if the number of the waiting drivers (riders) is very large, then more riders (drivers) should be sent here to avoid reneging penalties. 
The JLQ policy is proposed along this idea, which assigns a newly arrived traveler to the match with the largest number of travelers waiting to share a ride with the new traveler.
We refer to the baseline policy that myopically assigns travelers as the \emph{myopic} policy in the sequel. 
When the myopic policy is employed and a traveler of type $\ell\in{\color{black}\mL}$ arrives at time $t$, the traveler will be assigned to the match $e$ (incident to the vertex $\ell$) with the highest immediate expected return. The immediate expected return is considered as the difference between the immediate expected rewards gained by selecting and not selecting this match; that is,  if $\ell^D(e)=\ell$, then
the immediate expected return is set to $R({\color{black}N^{\text{myopic}}_e(t)},1,1)-R({\color{black}N^{\text{myopic}}_e(t)},0,1)$;
otherwise, $R({\color{black}N^{\text{myopic}}_e(t)},1,1)-R({\color{black}N^{\text{myopic}}_e(t)},1,0)$. 
Tie cases are broken by choosing the least label of the match.

\begin{figure}[t]
\centering
\begin{minipage}[]{0.33\textwidth}
\centering \vspace{-0.3cm}
\includegraphics[width=\linewidth]{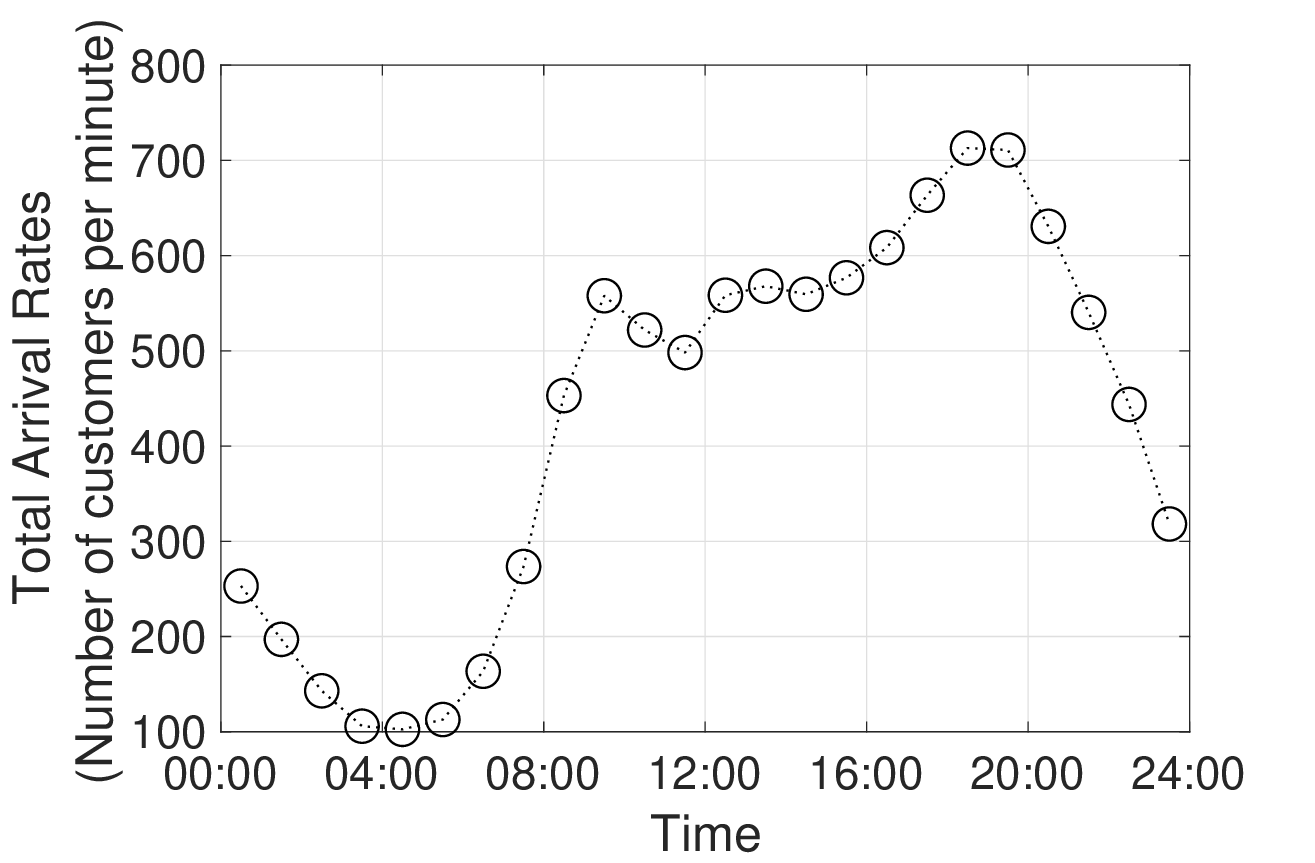}\vspace{0.15cm}
\caption{Arrival rates averaged in each hour.}
\label{fig:arrival_rates}
\end{minipage}
\begin{minipage}[]{0.66\textwidth}
\centering
\subfigure[]{\includegraphics[width=0.46\linewidth]{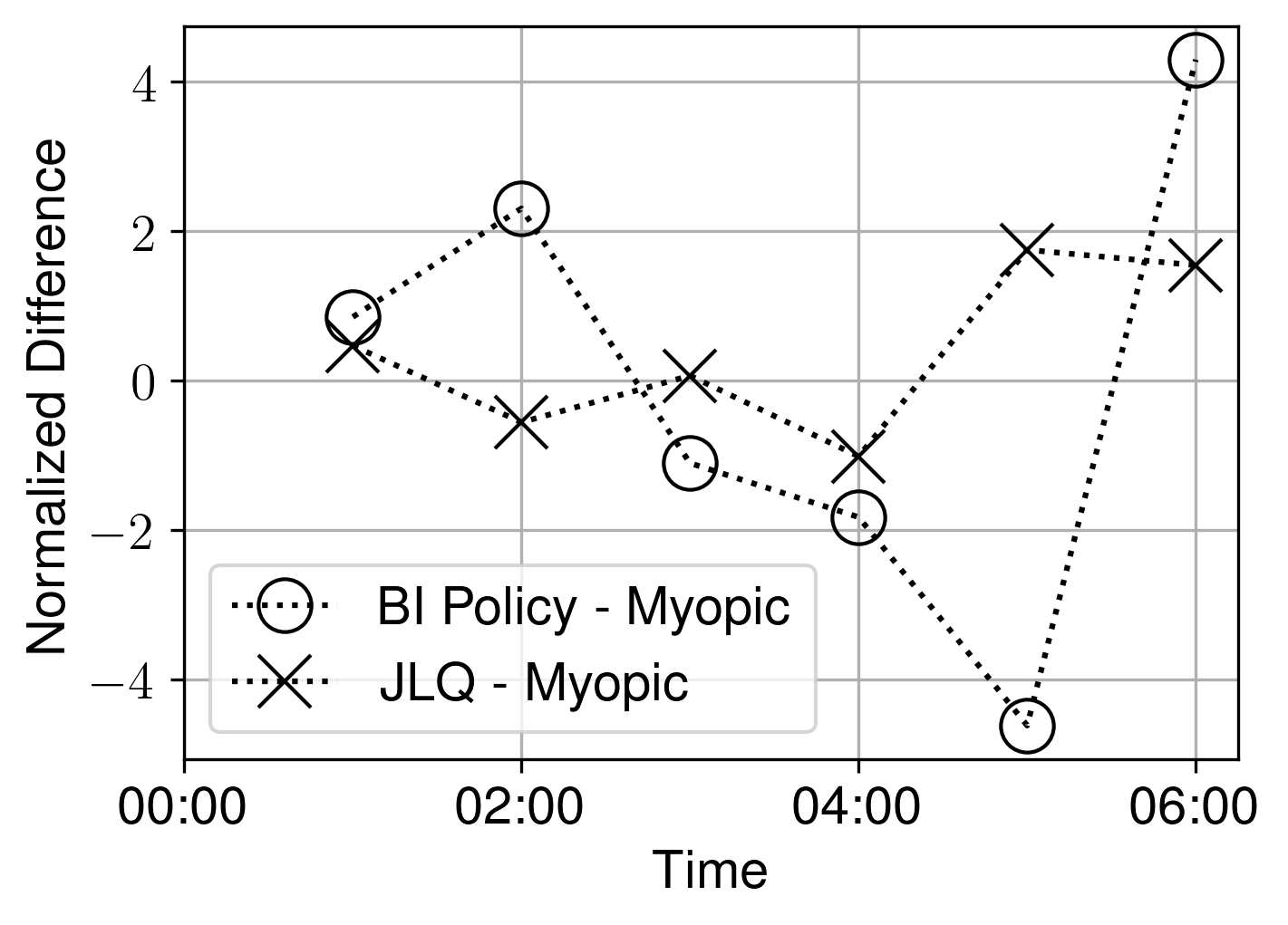}\label{fig:melb_eta4_diff:1}}
\subfigure[]{\includegraphics[width=0.49\linewidth]{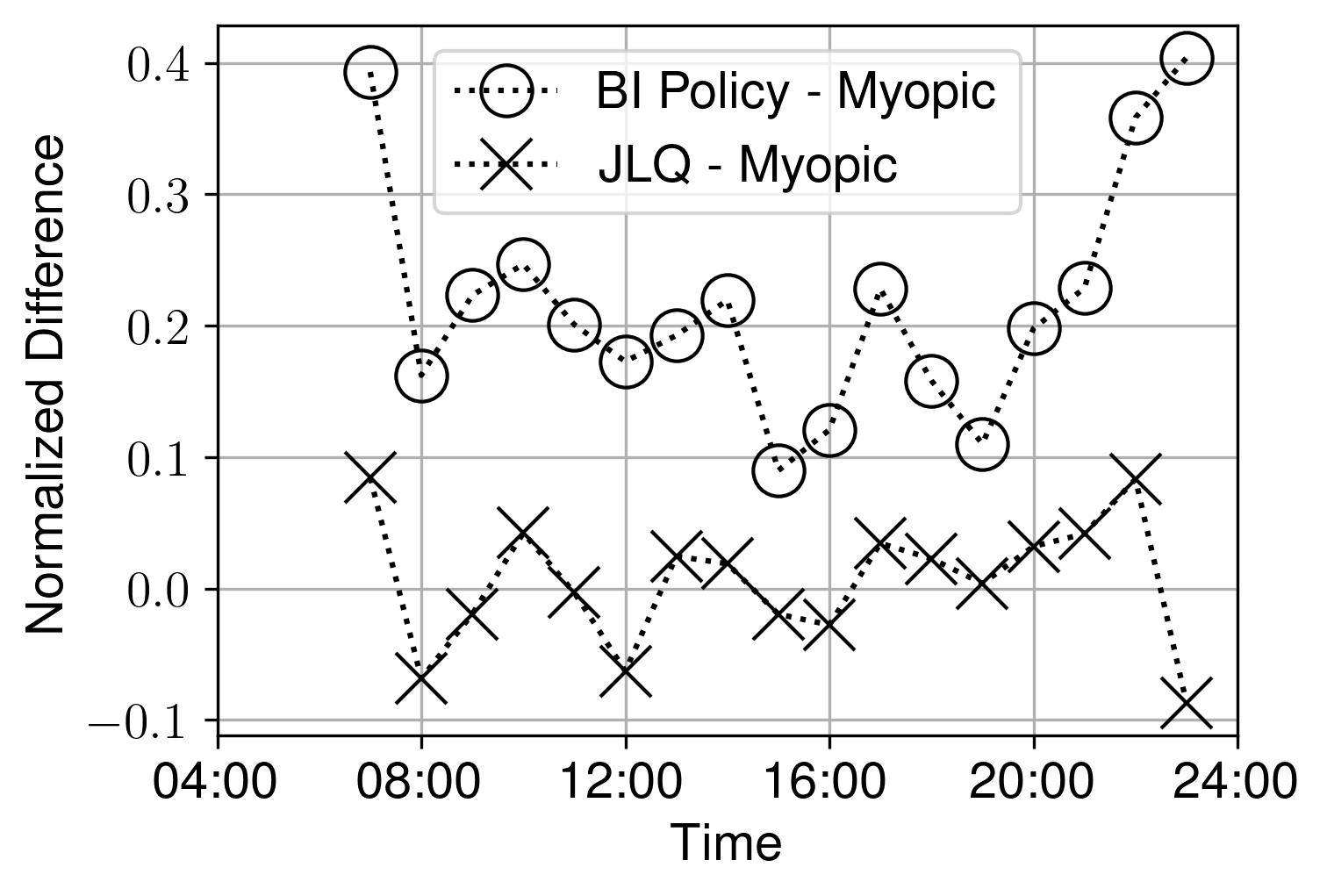}\label{fig:melb_eta4_diff:2}}
\vspace{-0.5cm}\caption{Relative difference of specific policies to JLQ with respect to average rewards for each time slot when $\zeta=4$: (a) from 0:00 to 6:00; and (b) from 6:00 to 24:00.}
\label{fig:melb_eta4_diff}
\end{minipage}
\vspace{-0.3cm}
\end{figure}

\begin{figure}[t]
\centering
\begin{minipage}[]{0.33\textwidth}
\centering\vspace{-0.2cm}
\includegraphics[width=\linewidth]{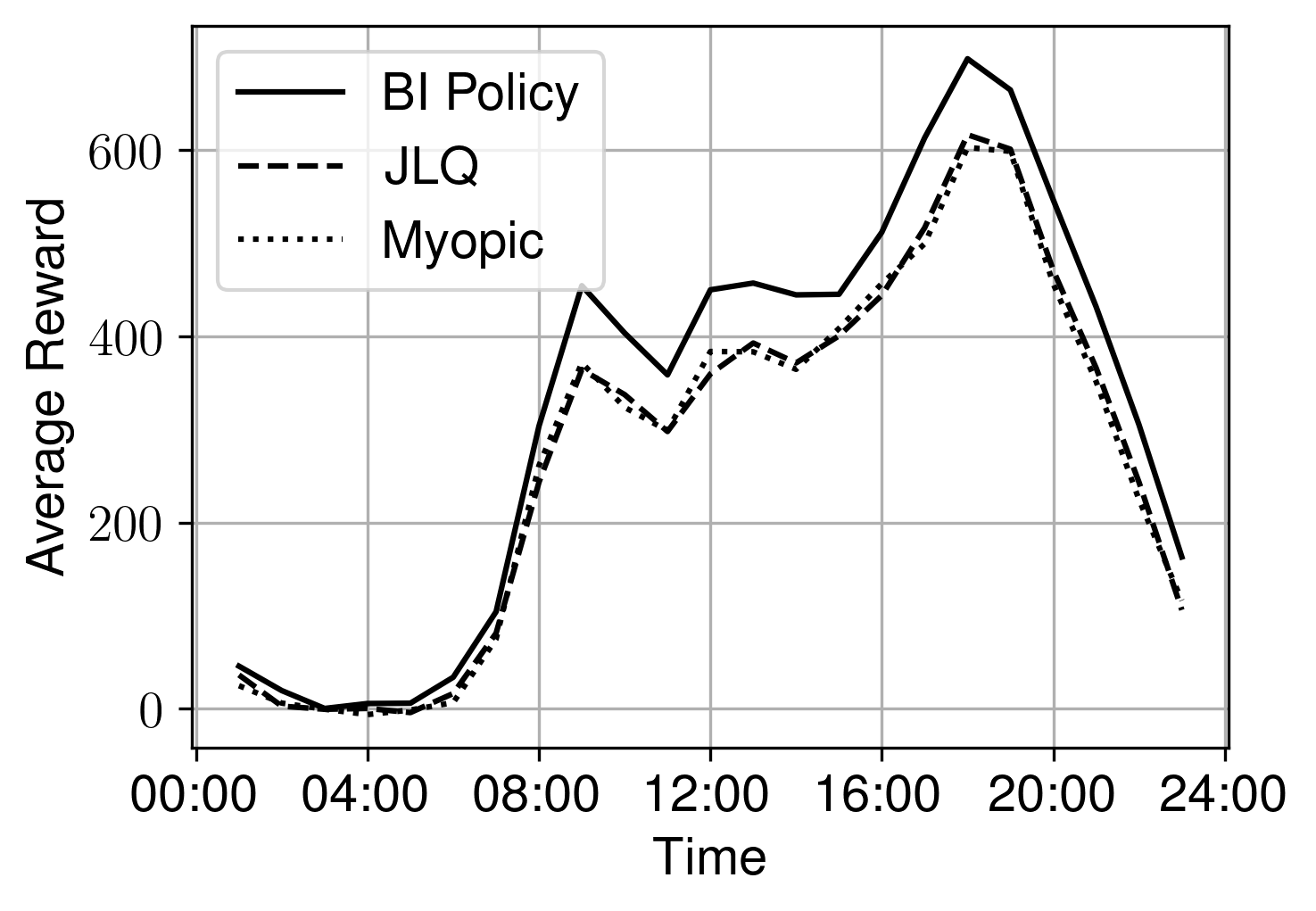}\vspace{0.1cm}
\caption{Reward rates under different policies against timeline when $\zeta=4$.}
\label{fig:melb_eta4}
\end{minipage}
\begin{minipage}[]{0.66\textwidth}
\subfigure[]{\includegraphics[width=0.49\linewidth]{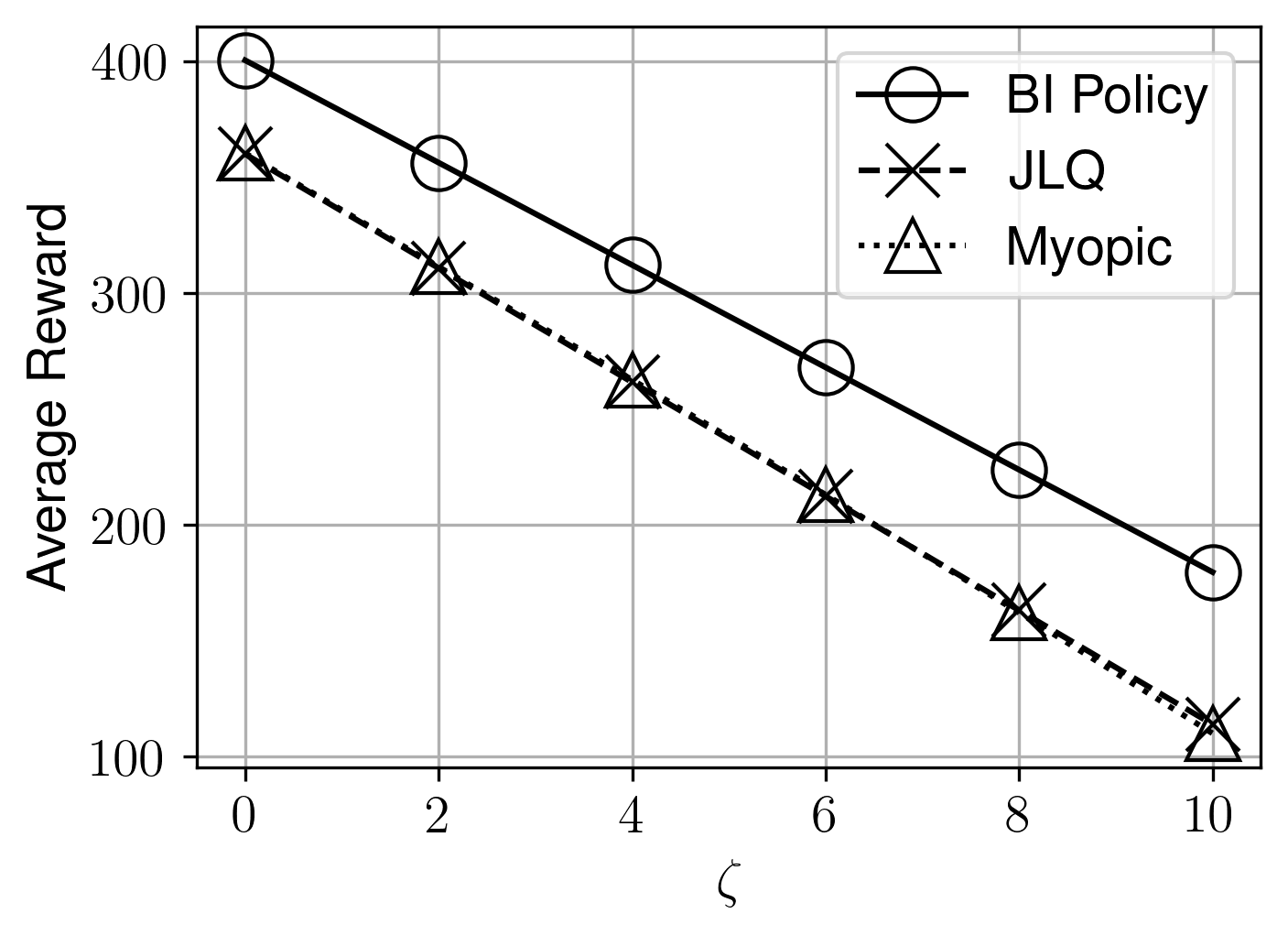}\label{fig:melb_reward}}
\subfigure[]{\includegraphics[width=0.48\linewidth]{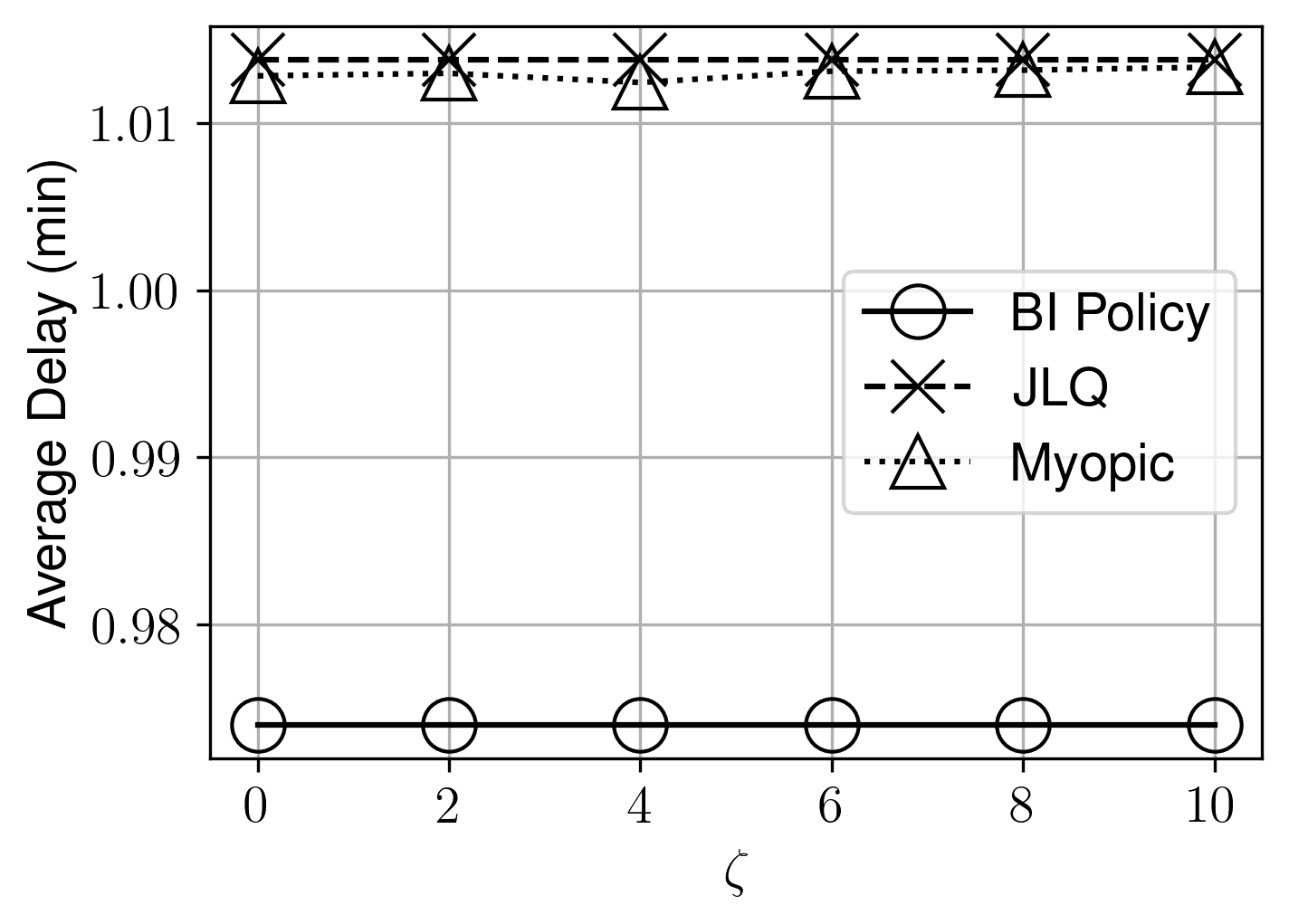}\label{fig:melb_delay}}
\vspace{-0.5cm}\caption{Average reward and delay for the Melbourne network gained under different policies against penalty parameter $\zeta$: (a) the relative difference of specified policy to the myopic with respect to the average reward; and (b) the average delay.}
\label{fig:melb_reward_delay}
\end{minipage}
\vspace{-0.5cm}
\end{figure}

We construct the \emph{Melbourne network} by dividing the Melbourne area (within the 5km radius from Melbourne CBD) into $20$ districts based on Local Government Areas (LGAs) defined in \cite{ABS2018}.
These LGAs are illustrated in Figure~\ref{fig:melbourne_example:demands}. 
In this section, the expected travel distance between addresses in different districts is approximated by the distance between the \emph{center points} of the corresponding districts. 
We classify drivers or riders according to their sources and destinations, and furthermore, we classify them by genders and smoking preferences. The sex ratio is set to be 1, while the ratio of smokers to non-smokers is $3/17\approx 0.176$ according to \citep{Smoke}. So, drivers or riders may specify if they would not ride with someone with a different gender or smoking preference. In total, we have $L=3040$ traveler types.
The reneging rates, reneging penalties,  and matching mechanism and rewards for two different traveler types are described in~\ref{app:uniform_network_settings} of the supplementary material. 
According to the geographic deployment of Melbourne destinations, in total $|\mE|=E=7400$ matches can be generated for the $L=3040$ traveler types.
We provide two examples of such matches in Figures~\ref{fig:melbourne_example:2}-\ref{fig:melbourne_example:3}.
The red and blue districts are the source and destination, respectively, of a driver while the light red and light blue districts are the sources and destinations, respectively, of the riders that can match this driver. 
Overall, travelers with neighboring sources and destinations are more willing to share the rides.

Consider the arrival rates for travelers of type $\ell$ to be time-dependent and positively affected by the \emph{weighted population} of their source and destination districts. 
More precisely, let $src(\ell)$ and $dst(\ell)$ represent the labels of the source and destination LGAs, respectively, of the traveler type $\ell$, and let
$\lambda_{\ell}= \frac{\Omega(\ell)}{\sum_{\ell'\in{\color{black}\mL}}\Omega(\ell')}\bar{\Lambda}(t)$,
with $\Omega(\ell)\coloneqq \sqrt{\omega_{src(\ell)}p_{src(\ell)}\omega_{dst(\ell)}p_{dst(\ell)}}$ for all $\ell\in{\color{black}\mL}$, 
where $p_{j}$ is the population in the LGA district labeled by $j$, $w_j\geq 1$ is a weighting factor for the district $j$ representing the attractiveness of a district to travelers, and $\bar{\Lambda}(t)$ is the total arrival rate of the ride-matching system at time $t$.
The detailed settings of $w_j$ and $p_j$ for all the districts are provided in~\ref{app:melbourne_case} of the supplementary material, for which the population data for each district used in this section is provided by \cite{ABS2018}. 
The resulting $w_j p_j$, also referred to as the weighted population, of all the LGAs are illustrated in Figure~\ref{fig:melbourne_example:demands} with approximated border lines.
Based on the data-set of the taxi trips for May 2016, reported to the City of Chicago \citep{Chicago} (due to lack of local data in Melbourne),  the total arrival rate $\bar{\Lambda}(t)$, estimated as the number of arrived travelers per minute, averaged in each hour is presented in Figure~\ref{fig:arrival_rates}.


In Figures~\ref{fig:melb_eta4_diff} and~\ref{fig:melb_eta4}, we demonstrate the average rewards for BI, JLQ and the myopic policy and the relative difference of BI and  JLQ to the myopic policy with respect to the average reward per each hour, respectively, for the Melbourne network. 
In particular, the relative difference of policy
$\phi_1$ to policy $\phi_2$ is $(\rho^{\phi_1}-\rho^{\phi_2})/|\rho^{\phi_2}|$, where $\rho^{\phi}$ represents the average reward (per hour) under policy $\phi$.
The BI policy significantly outperforms the other policies for almost all the tested time slots. 
In particular, for idle hours between 00:00 and 08:00, the differences between the tested policies are not significant; while, during relatively busy hours between 08:00 and 24:00, BI policy achieves clearly higher reward rates, which can be up to $40\%$ higher than that of JLQ and the myopic policy. We remark that the relative differences between BI and the myopic policy at some time slots are negative in Figure~\ref{fig:melb_eta4_diff:1} because the myopic produces negative rewards, and in fact, BI always outperforms the myopic policy.
With the same settings as the simulations presented in Figure~\ref{fig:melb_eta4_diff} and Figure~\ref{fig:melb_eta4}, in Figure~\ref{fig:melb_reward}, we present the reward averaged for the entire $24$ hours and the relative difference of the average reward against $\zeta$. 
Note that the average reward for Figure~\ref{fig:melb_reward} is the ratio of the reward accumulated within the $24$ hours to the time horizon (that is, $24$ hours) of the tested process.
The BI policy achieves clearly higher rewards than the myopic policy and JLQ, and the relative difference of BI to either JLQ or the myopic policy is at least $12\%$ with respect to the average reward for all the tested $\zeta$.
In Figure~\ref{fig:melb_delay}, given the significantly higher average reward achieved by BI, we consider the average delay (waiting time) of travelers and observe compatible average delay under the BI policy, when compared to the other two policies.

\vspace{-0.7cm}
\section{Conclusion}\label{sec:conclusion}
\vspace{-0.5cm}
We have formulated the multi-rider matching problem as an MDP consisting of parallel sub-processes, where travelers dynamically arrive with potentially different travel preferences, sources and destinations.
We have considered the reneging behaviors of the travelers without assuming any maximal waiting time. 
This problem, referred to as the RMP, is complicated by the heterogeneous travelers who may share rides with each other, the reneging behaviors and the large problem size.

We have proposed the BI policy by extending the Whittle relaxation technique to the RMP. 
The storage and computational complexity of implementing the BI policy are linear and logarithmic, respectively, to the number of candidate matches between different riders and drivers.
The BI policy prioritizes travelers and the candidate matches according to a ranking of their bivariate indices, which intuitively represent the marginal profits of selecting certain riders and drivers to match each other.
We have proved that, in a special case, the bivariate indices lead to a policy optimal to the relaxed version of the RMP.
For the general case, the effectiveness of the BI policy has been demonstrated through extensive numerical simulations in real-world settings. Compared to the baseline policies, BI achieves significantly higher average rewards with a compatible average waiting time of the served travelers.
This is also the first attempt at adopting and extending the Whittle relaxation technique to the problem, where sub-processes with multiple action variables subject to different action constraints are inevitable.

{\color{black}Building on the insights gained in this paper, delving into specific application scenarios with practical configurations in greater detail would be beneficial for future research.
For instance, this paper considers a given rule on the eligibility of matching a driver and a rider (the given bipartite graph $G$), which plays an important role in approximating optimality. 
The matching eligibility integrates travellers' diverse preferences, traffic conditions, routing topology, and all the other complex environmental factors. In-depth exploration of more practical and flexible matching mechanisms will enhance the model’s interpretability and applicability in real-world contexts.
Moreover, when we lack an accurate estimation on the matching mechanism, the mis-specified matching eligibility may adversely impact the performance of the BI policy.
Additional studies with only partial knowledge of the bipartite graph $G$ (matching eligibility) and/or other system parameters will lead to more flexible and practical results, thereby deepening the practical implications of all related studies.

On the other hand, for computing the bivariate indices, for the general case where the hypothesis of Proposition~\ref{prop:priority} is unlikely to satisfy, further studies can explore appropriate estimations of $\eta_2\in\mathbb{R}$, based on which we can solve the $(\bm{\eta},e)$-process and obtain $\Gamma^{n,e}(\eta_2)$ for all $\Gamma^{n,e}\in\mathscr{B}^1(n,e)$, $n\in\tilde{\mathscr{N}}\backslash\{N\}$ and $e\in\mE$. 
These $\Gamma^{n,e}(\eta_2)$ can serve as adapted bivariate indices for the BI policy that further improves the overall performance. }

\vspace{-0.7cm}
\section{Acknowledgement}
\vspace{-0.5cm}
We would like to acknowledge Dr. Kai Huang for his help in this work.

\vspace{-0.7cm}

\bibliographystyle{elsarticle-harv} 
\bibliography{references/IEEEabrv,references/bib1,references/rideshare}

\clearpage
\setcounter{page}{1}

\appendix
\section{Example of the Matching Mechanism}\label{app:example:2}

\begin{figure}[t]
\centering
\begin{minipage}[]{0.49\textwidth}
\centering
\includegraphics[width=0.8\linewidth]{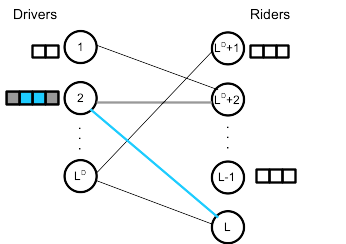}
\vspace{-0.3cm}\caption{An example of the bipartite graph for the matches between drivers and riders.\label{fig:bipartite}
}
\end{minipage}
\begin{minipage}[]{0.49\textwidth}
\centering
\includegraphics[width=0.55\linewidth]{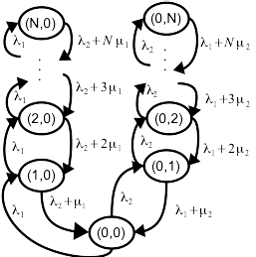}
\vspace{-0.2cm}\caption{The underlying Markov chain of the processes $\{(N_{\ell^D(e),e}(t),N_{\ell^R(e),e}(t)), t\geq 0\}$ associated with the match $e$, where $\lambda_1=\lambda_{\ell^D(e)}a_{\ell^D(e),e}({\color{black}n}(t))$, $\lambda_2 = \lambda_{\ell^R(e)}a_{\ell^R(e),e}({\color{black}n}(t))$, $\mu_1 = \mu_{\ell^D(e),e}$ and $\mu_2 = \mu_{\ell^R(e),e}$. 
\label{fig:markov_chain}
}
\end{minipage}
\vspace{-0.5cm}
\end{figure}

\begin{example}\label{example:2}
Consider a simple example in Figure~\ref{fig:bipartite}, the left and right vertices represent driver and rider types, respectively, and the edges indicate the eligibility of matching the incident drivers and riders. 
In this figure, the edges are labeled by $e=1,2,\ldots,5$.
In this example, the drivers of type $1$ can only match riders of type $2$. 
When there is no rider of type $2$, all the two arrived drivers of type $1$ have to wait in the queue ($N_{1,1}(t)=2$).
The drivers of type $2$ can be matched to either riders of types $L^D+2$ (grey edge) or $L$ (blue edge).
Based on an employed scheduling policy, the first and forth drivers (grey squares) of type $2$ are waiting for riders of type $L^D+2$, and the other two (blue squares) awaiting riders of type $L$.
At this moment, $N_{2,2}(t) =2$ and $N_{2,3}=2$.

On the other hand, the riders of type $L^D+1$ have to wait in the queue until a driver of type $L^D$ arrives ($N_{L^D+1,4}=3$). 
If there is more than one eligible match (edge) for a rider type $\ell$, then, similar to drivers of type $2$, different riders of the type $\ell$ may queue there awaiting drivers of different types.
\end{example} 

For $e\in{\color{black}\mE}$, the state transition rates of the underlying process $\bigl\{{\color{black}N_e(t)}, t\geq 0\bigr\}$ are presented in Figure~\ref{fig:markov_chain}.

{\color{black}\section{Computing Whittle Indices}\label{app:whittle}
For the $e$-process, the pseudo-code for computing Whittle indices, with respect to the driver, is provided in Algorithm~\ref{algo:whittle}.
Recall that the Whittle indices can be computed/approximated through a couple of conventional techniques, such as the bisection method in \citep{fu2019towards,fu2021large} and the marginal productivity (MP) index in \citep{nino2001restless,nino2002dynamic,nino2007dynamic,nino2020verification}, or exists in a closed form for a range of special cases in \citep{fu2023restless}.
Here, we use the bisectioning method used in \cite{fu2021large}.
The computational complexity of Algorithm~\ref{algo:whittle} is $O(\lvert\tilde{\mathscr{N}}\rvert^2 C_{\text{VI}}(\sigma)\log_{\sigma}M)$, where $C_{\text{VI}}(\sigma)$ is the number of iterations for the value functions $V_e(n)$ ($n\in\tilde{\mathscr{N}}\backslash\{N\}$ in Lines~\ref{line:value_iteration:1}-\ref{line:value_iteration:2}, and $\log_{\sigma}M$ is the complexity for bisectioning $\eta^1_e(n)$, $n\in\tilde{\mathscr{N}}\backslash\{N\}$.
The iteration method is classic for approximating value functions of MDP and is in general efficient~\cite{ross1992applied}.
Recall that the indices $\eta^1_e(n)$, $(n,e)\in\Bigl(\tilde{\mathscr{N}}\backslash\{N\}\Bigr)\times\mE$, are pre-computed offline, and can be computed in parallel for all $e\in \mE$.
Also, since the indices are computed independently per each $e\in\mE$, for a system with newly added eligible matches, the computation of the indices for the new matches does not demand re-calculation of the original ones.

The computation of the indices for the rider, $\eta^2_e(n)$, is completed along the same lines of Algorithm~\ref{algo:whittle} except swapping the second and third arguments of functions $R_e$ and $u_e$ in Lines~\ref{line:swap:1}, \ref{line:swap:2}, and \ref{line:swap:3}. 
It achieves the same computational complexity as that of $\eta^1_e(n)$.

With given indices $\eta^1_e(n)$ and $\eta^2_e(n)$ for all $(n,e)\in \Bigl(\tilde{\mathscr{N}}\backslash\{N\}\Bigr)\times \mE$ and $\Bigl(\tilde{\mathscr{N}}\backslash\{-N\}\Bigr)\times \mE$, respectively, the BI policy is implemented by searching for the eligible match $e$ among $\mE_{\ell}$ upon a new arrival of traveller of type $\ell\in\mL$.
Here, we consider a max-heap-based method to rank the eligible matches for BI.
For each $\ell\in\mL^D$ (or $\ell\in\mL\backslash\mL^D$), we maintain a max heap $\mathcal{H}_{\ell}(t)$ of all the eligible matches $e\in\mE_{\ell}$ with respect to their bivariate indices $\eta^1_e\bigl(N^{\text{BI}}_e(t)\bigr)$ (or $\eta^2_e\bigl(N^{\text{BI}}_e(t)\bigr)$).
Let $\upsilon_{\ell}(t)$ indicate the match $e$ with the highest bivariate index in $\mE_{\ell}$; that is, $\upsilon_{\ell}(t)$ is equal to the root node of the max heap $\mathcal{H}_{\ell}(t)$.
If all the matches $e\in\mE_{\ell}$ have $N$ travellers of type $\ell$ waiting, then $\upsilon_{\ell}(t) = -1$.
Upon an arrival of a type-$\ell$ traveller, it will be assigned to the match $\upsilon_{\ell}(t)$, or, if $\upsilon_{\ell}(t)=-1$, this traveller will be rejected.
Let $\pmb{\upsilon}(t)\coloneqq (\upsilon_{\ell}(t):\ell\in\mL)$.
If the traveller is assigned to a match $e\in\mE_{\ell}$, then the state variable $N^{\text{BI}}_e(t)$ changes, incurring updates of $\mathcal{H}_{\ell}(t)$ and $\mathcal{H}_{\ell'}(t)$ where $\ell'$ is the other vertex in $G$ incident edge $e$.
The vector $\pmb{\upsilon}(t)$ changes accordingly.
We provide pseudo-code in Algorithm~\ref{algo:arrival} for updating the max heaps and $\pmb{\upsilon}(t)$ upon a state transition $N^{\text{BI}}_{\upsilon_{\ell}(t)}(t)$ that is caused by an arrival of type-$\ell$ traveller.
The steps in Algorithm~\ref{algo:arrival} has computational complexity $O(\log |\mE_{\ell}| + \log|\mE_{\ell'}|)$ (updating the max heaps) for each $\ell\in\mL$, where $\ell'$ is the other vertex incident edge $\upsilon_{\ell}(t)$, and storage complexity $O(|\mE_{\ell} + \mE_{\ell'}|)$.
}

\IncMargin{1em}
\begin{algorithm*}\small
\linespread{0.5}\selectfont
\color{black}
\SetKwFunction{WhittleIndices}{$(\bm{\eta}_e)\gets$ WhittleIndices}
\SetKwProg{Fn}{Function}{:}{End}
\SetKwInOut{Input}{Input}\SetKwInOut{Output}{Output}
\SetAlgoLined
\DontPrintSemicolon
\Input{$e\in \mE$, $\beta\in [0,1)$ and $\sigma \in\mathbb{R}_+$} 
\tcc*{$\beta\uparrow 1$ is the hyper-parameter for Blackwell optimality.}
\tcc*{$\sigma\downarrow 0$ is the hyper-parameter used to indicate the accuracy.}
\Output{$\bm{\eta}_e\in\mathbb{R}^{|\tilde{\mathcal{N}}|}$}

\Fn{\WhittleIndices{$e,\beta,\sigma$}}{
	\For{$n\in\tilde{\mathcal{N}}\backslash\{N\}$}{
	 	$\eta_0,\eta_1\gets -M, +M$ \tcc*{Lower and upper bounds for the Whittle index of state $\bm{n}$, where $M\in\mathbb{R}_0$ is sufficiently large.}
    	\While{1}{   					
		    $V_e(n')\gets 0 $ for all $n'\in\tilde{\mathcal{N}}$\; 
			$\eta\gets (\eta_0+\eta_1)/2$ \tcc*{Bisectioning}
			\While{1\label{line:value_iteration:1}}{ 
				$\hat{V}_e(n')\gets \max\limits_{a\in\{0,1\}}\left\{\frac{R_e(n',a,1)-a\eta }{u_e(n',(a,1))} 
+\beta\sum\limits_{n''\in \tilde{\mathcal{N}}}p_e^a(n',n'')V_e(n'')\right\}$ for all $n'\in\tilde{\mathcal{N}}\backslash\{n^*\}$\label{line:swap:1}\;
				\If{$|\hat{V}_e(n')-V_e(n)| < \sigma$ for all $n'\in\tilde{\mathcal{N}}\backslash\{n^*\}$}{
					{\bf Break}\;
				}
				$V_e(n')\gets \hat{V}_e(n')$ for all $n'\in\tilde{\mathcal{N}}\backslash \{n^*\}$ \tcc*{Value iteration}
			\label{line:value_iteration:2}}
			$V_1 \gets \frac{R_e(n,1,1)-\eta }{u_2(n,(1,1))} 
+\beta\sum\limits_{n'\in \tilde{\mathcal{N}}}p_e^1(n,n')V_e(n')$\label{line:swap:2}\;
			$V_0\gets \frac{R_e(n,0,1)}{u_e(n,(0,1))} 
+\beta\sum\limits_{n'\in \tilde{\mathcal{N}}}p_e^0(n,n')V_e(n')$\label{line:swap:3}\;
			\uIf{$|V_1-V_0|< \sigma$}{
				{\bf Break}  
									\tcc*{The accuracy is guaranteed by hyper-parameter $\sigma$.}
			}\uElseIf{$V_1 > V_0$}{
				$\eta_0 \gets \eta$ \;
			}\Else{
				$\eta_1 \gets \eta$ \;
			}
		}
		$\eta^1_e(n)\gets \eta$\;
    }    

	$\eta^1_e(n^*) \gets 0$\;
	\Return $(\eta^1_e(n):~n\in\tilde{\mathcal{N}})$\;
}
\caption{Pseudo-code for computing Whittle indices to the $e$-process with respect to the driver.}\label{algo:whittle}
\end{algorithm*}
\DecMargin{1em}

\IncMargin{1em}
\begin{algorithm*}\small
\linespread{0.5}\selectfont
\color{black}
\SetKwFunction{UpdatingUponArrival}{}
\SetKwProg{Fn}{Function}{:}{End}
\SetKwInOut{Input}{Input}\SetKwInOut{Output}{Output}
\SetAlgoLined
\DontPrintSemicolon
\Input{$\ell\in \mL$ and the updated state vector $\bm{N}^{\text{BI}}(t)$.} 
\Output{$\pmb{\upsilon}(t)$, $\mathcal{H}_{\ell}$ and $\mathcal{H}_{\ell'}$ where $\ell'$ is the other vertex in $G$ incident edge $e=\upsilon_{\ell}(t)$.}
\Fn{\UpdatingUponArrival{$\ell,\bm{N}^{\text{BI}}(t)$}}{
    Let $\ell'\neq\ell$ represent the vertex in $G$ that incidents edge $e=\upsilon_{\ell}(t)$.\;
    $e\gets \upsilon_{\ell}(t)$\;
    \uIf{$\ell\in \mL^D $ and $N^{\text{BI}}_e(t) = N$}{
        Remove $e$ from max heap $\mathcal{H}_{\ell}$\;
    }\uElseIf{$\ell\in \mL\backslash\mL^D$  and    $N^{\text{BI}}_e(t) = -N$}{
        Remove $e$ from max heap $\mathcal{H}_{\ell}$\;
    }\Else{
        Update $\mathcal{H}_{\ell}$ according to the updated $N^{\text{BI}}_e(t)$\;
    }
    $\upsilon_{\ell}(t)\gets$ the root node of max heap $\mathcal{H}_{\ell}(t)$\;

    \tcc*{The max heap for $\ell'$ is updated in the same way.}\;

    \uIf{$\ell'\in \mL^D $ and $N^{\text{BI}}_e(t) = N$}{
        Remove $e$ from max heap $\mathcal{H}_{\ell'}$\;
    }\uElseIf{$\ell'\in \mL\backslash\mL^D$  and    $N^{\text{BI}}_e(t) = -N$}{
        Remove $e$ from max heap $\mathcal{H}_{\ell'}$\;
    }\Else{
        Update $\mathcal{H}_{\ell'}$ according to the updated $N^{\text{BI}}_e(t)$\;
    }
    $\upsilon_{\ell'}(t)\gets$ the root node of max heap $\mathcal{H}_{\ell'}(t)$\;
    \Return $\pmb{\upsilon}(t)$, $\mathcal{H}_{\ell}(t)$, and $\mathcal{H}_{\ell'}(t)$;
}
\caption{Updating $\pmb{\upsilon}(t)$ and the max heaps.}\label{algo:arrival}
\end{algorithm*}
\DecMargin{1em}

\section{Settings for Simulations in Figure~\ref{fig:example_1} }\label{app:settings:example_1}

Consider a simple system with $L=2$ traveler types and $|\mE|=1$ match. 
The only match $e$ satisfies $\ell^D(e) = 1$ and $\ell^R(e)=2$.
The arrival rates, $\lambda_{\ell}$,  of all traveler types $\ell\in{\color{black}\mL}$ are set to be $5$; and the reneging rates $\mu_{\ell,e}$ ($\ell\in{\color{black}\mL},e\in{\color{black}\mE}$) are set to be $1$. The upper bound of traveler number of each type waiting for the match, $N$, is considered to be $5$ so that the arrival rates of new travelers and the total reneging rates of waiting travelers are compatible.
Also, the reward earned by serving a couple of travelers through the match $e$ is $R_e=10$. In this example, we do not consider the reneging penalty; that is, $C_{\ell,e}=0$ ($\ell\in{\color{black}\mL},e\in{\color{black}\mE}$).
The curves presented in Figure~\ref{fig:example_1} are the simulated results of the sub-process associated with the only match.

\section{Proof of Proposition~\ref{prop:priority}}\label{app:priority}

\begin{lemma}\label{lemma:priority}
For $e\in{\color{black}\mE}$ and ${\color{black}n}\in\tilde{\mathscr{N}}\backslash\{{\color{black}N}\}$, if $\mu_{\ell^D(e)}=\mu_{\ell^R(e)} = 0$, then
there exists $E\in\mathbb{R}$ and $\Gamma\in\mathcal{B}_1({\color{black}n},e)$ such that, for any $\eta_2 < E$, 
\begin{equation}\label{eqn:lemma:priority}
\Gamma(\eta_2) = 
    -\frac{\lambda_{\ell^D(e)}}{\lambda_{\lambda^R(e)}}(\eta_2 + g_e(\bm{\eta})) +\lambda_{\ell^D(e)}R_e,  
\end{equation}
where $g_e(\bm{\eta})\in\mathbb{R}$ is equal to the maximal average reward of the $(\bm{\eta},e)$-process.
\end{lemma}
\proof{Proof.}
For $e\in{\color{black}\mE}$ and any given $\bm{\eta}\in\mathbb{R}^2$, the action variables $\bm{a}$ that maximizes right hand side of equation \eqref{eqn:bellman_duplex} for all ${\color{black}n}\in\tilde{\mathscr{N}}\backslash\{{\color{black}0}\}$ are the same as
\begin{equation}\label{eqn:bellman_1}
\max\limits_{\bm{a}\in\mathcal{A}({\color{black}n})}\left\{
\frac{R_e({\color{black}n},a_1,a_2)-a_1\eta_1-a_2\eta_2 -g_e(\bm{\eta})}{u_e({\color{black}n},\bm{a})} 
+\sum\limits_{\bm{n'}\in \tilde{\mathscr{N}}}p_e^{\bm{a}}({\color{black}n},{\color{black}n}')V^{\bm{\eta}}_e({\color{black}n}')
\right\},
\end{equation}
where  $g_e(\bm{\eta})\in\mathbb{R}$  is a given real number equal to the maximized average reward of the $(\bm{\eta},e)$-process.
Since $R_e({\color{black}n},a_1,a_2)$ is bounded and $u_e({\color{black}n},\bm{a})>0$ for all ${\color{black}n}\in \tilde{\mathscr{N}}\backslash\{{\color{black}0}\}$, for given $\eta_1\in\mathbb{R}$, there exists $E\in\mathbb{R}$ such that, for all $\eta_2<E$, the optimal vector $\bm{a}$ which maximizes right hand side of \eqref{eqn:bellman_1} satisfies $a_2=a_2(n_2)$, where \vspace{-0.5cm}
\begin{equation}\label{eqn:fixed_a2}
a_2({\color{black}n}) \coloneqq  \begin{cases}
1,&\text{if } {\color{black}n > -N}, \\
0,&\text{otherwise}. 
\end{cases}
\end{equation}
\vspace{-0.5cm}
That is, for all $\eta_2<E$, the maximization in \eqref{eqn:bellman_1} becomes
\begin{equation}\label{eqn:bellman_2}
\max\limits_{a_1\in\bigl[0,\mathds{1}\{{\color{black}n}<N\}\bigr]}\biggl\{
\frac{R_e({\color{black}n},a_1,a_2(n))-g_e(\bm{\eta})-a_2({\color{black}n})\eta_2-a_1\eta_1}{u_e({\color{black}n},(a_1,a_2({\color{black}n})))} 
+\sum\limits_{{\color{black}n'}\in \tilde{\mathscr{N}}}p_e^{a_1}({\color{black}n},{\color{black}n}')V^{\bm{\eta}}_e({\color{black}n}')
\biggr\},
\end{equation}
where $p^{a_1}_e({\color{black}n},{\color{black}n}')$ is defined to be equal to $p^{\bm{a}}_e({\color{black}n},{\color{black}n}')$ with $a_2={\color{black}a_2(n)}$.
Along similar lines, we can show that, for any $\eta_1\in\mathbb{R}$, there exists $E\in\mathbb{R}$ such that, for all $\eta_2<E$, a vector $\bm{a}=(a_1,a_2)$ maximizes the right-hand side of \eqref{eqn:bellman_zero} if and only if
$a_2=a_2({\color{black}n})$ and $a_1$ maximizes the right hand side of
\begin{equation}\label{eqn:bellman_3}
H^{\bm{\eta}}_e=\max\limits_{a_1\in[0,1]}\biggl\{
\frac{R_e({\color{black}0},a_1,1)-g_e(\bm{\eta})-\eta_2-a_1\eta_1}{u_e({\color{black}0},(a_1,1))} 
+\sum\limits_{\bm{n'}\in \tilde{\mathscr{N}}}p_e^{a_1}({\color{black}0},{\color{black}n}')V^{\bm{\eta}}_e({\color{black}n}')
\biggr\}.
\end{equation}
The $H^{\bm{\eta}}_e$ is equal to zero if and only if $g_e(\bm{\eta})$ is equal to the maximized long-run average reward of the $(\bm{\eta},e)$-process. 

It remains to show that, if $\mu_{\ell^D(e)} = \mu_{\ell^R(e)} = 0$, then there exists $\Gamma\in\mathscr{B}_1({\color{black}n},e)$ such that, for all $\eta_2<E$, \eqref{eqn:lemma:priority} is satisfied.
We can prove it by considering four cases: ${\color{black}n}\in\mathscr{N}_-\coloneqq \{{\color{black}-N+1,-N+2,\ldots,-1}\}$, $\mathscr{N}_+\coloneqq \{{\color{black}1,2,\ldots,N-1}\}$, $\mathscr{N}_0\coloneqq \{{\color{black}0}\}$, and $\mathscr{N}_N\coloneqq \{{\color{black}-N}\}$.

In this proof, the following discussions all assume that $\mu_{\ell^D(e)} = \mu_{\ell^R(e)} = 0$.
In this context, if $\alpha^{\phi^*}_{\ell^D(e),e}({\color{black}n})=0$ then all the states ${\color{black}n'}$ with ${\color{black}n'>n}$ become transient, and it makes no difference when setting $\alpha^{\phi^*}_{\ell^D(e),e}({\color{black}n}')=0$ or $1$.
Let $a_1^*({\color{black}n})$, ${\color{black}n}\in\tilde{\mathscr{N}}\backslash\{{\color{black}N}\}$ represent the value of $\alpha^{\phi^*}_{\ell^D(e),e}({\color{black}n})$ under an optimal policy $\phi^*$ for the $(\bm{\eta},e)$-process. In particular, if $a_1^*({\color{black}n}) = 0$, then $a_1^*({\color{black}n}') = 0$ for all the states ${\color{black}n'}$ with ${\color{black}n'>n}$.

For ${\color{black}-n}\in\mathscr{N}_-$, $\eta_1=\Gamma(\eta_2)$ satisfies
\begin{multline}\label{eqn:lemma:priority:1}
    \frac{-g_e(\bm{\eta}) -\eta_2}{\lambda_{\ell^R(e)}} + V_e^{\bm{\eta}}({\color{black}-(n+1)})= \frac{\lambda_{\ell^D(e)}R_e - g_e(\bm{\eta})-\eta_2-\eta_1}{\lambda_{\ell^D(e)} + \lambda_{\ell^R(e)}} \\+ \frac{\lambda_{\ell^D(e)}}{\lambda_{\ell^D(e)} + \lambda_{\ell^R(e)}}V_e^{\bm{\eta}}({\color{black}-(n-1)}) +
    \frac{\lambda_{\ell^R(e)}}{\lambda_{\ell^D(e)} + \lambda_{\ell^R(e)}}V_e^{\bm{\eta}}({\color{black}-(n+1)}).
\end{multline}
The equation is equivalent to
\begin{equation}\label{eqn:lemma:priority:2}
    \eta_1 = -\lambda_{\ell^D(e)}\Bigl(V^{\bm{\eta}}_e({\color{black}-(n+1)})-V^{\bm{\eta}}_e({\color{black}-(n-1)})\Bigr)+\frac{\lambda_{\ell^D(e)}}{\lambda_{\ell^R(e)}}(g_e(\bm{\eta}) +\eta_2) +\lambda_{\ell^D(e)}R_e.
\end{equation}
From \eqref{eqn:bellman_2} and \eqref{eqn:lemma:priority:1}, we obtain 
\begin{equation}\label{eqn:lemma:priority:3}
    V^{\bm{\eta}}_e({\color{black}-n}) = \frac{-g_e(\bm{\eta}) -\eta_2}{\lambda_{\ell^R(e)}} + V_e^{\bm{\eta}}({\color{black}-(n+1)}).
\end{equation}
By plugging \eqref{eqn:lemma:priority:3} in \eqref{eqn:lemma:priority:2}, it follows that, for ${\color{black}-n}\in\mathscr{N}_-$, 
\begin{equation}\label{eqn:lemma:priority:4}
    \eta_1 = -\lambda_{\ell^D(e)}\Bigl(V^{\bm{\eta}}_e({\color{black}-n})-V^{\bm{\eta}}_e({\color{black}-(n-1)})\Bigr) + \lambda_{\ell^D(e)}R_2.
\end{equation}

For ${\color{black}-n}\in\mathscr{N}_-$, if $a_1^*({\color{black}-n}) = 0$,  then, based on the definition, $a_1^*({\color{black}-(n-1)})=0$. Together with \eqref{eqn:bellman_2}, it follows that
\begin{equation}\label{eqn:lemma:priority:5}
    V^{\bm{\eta}}_e({\color{black}-n})-V^{\bm{\eta}}_e({\color{black}-(n-1)}) = \frac{g_e(\bm{\eta})+\eta_2}{2}.
\end{equation}
Plugging in \eqref{eqn:lemma:priority:4}, we obtain \eqref{eqn:lemma:priority}.

Now we show that \eqref{eqn:lemma:priority} also holds for ${\color{black}n}\in\mathscr{N}_+$.
Similar to \eqref{eqn:lemma:priority:1}, for $\eta_1 = \Gamma(\eta_2)$, we have
\begin{multline}\label{eqn:lemma:priority:7}
    \frac{\lambda_{\ell^R(e)}R_e-g_e(\bm{\eta}) -\eta_2}{\lambda_{\ell^R(e)}} + V_e^{\bm{\eta}}(n-1)= \frac{\lambda_{\ell^R(e)}R_e - g_e(\bm{\eta})-\eta_2-\eta_1}{\lambda_{\ell^D(e)} + \lambda_{\ell^R(e)}} \\+ \frac{\lambda_{\ell^D(e)}}{\lambda_{\ell^D(e)} + \lambda_{\ell^R(e)}}V_e^{\bm{\eta}}(n+1) +
    \frac{\lambda_{\ell^R(e)}}{\lambda_{\ell^D(e)} + \lambda_{\ell^R(e)}}V_e^{\bm{\eta}}(n-1),
\end{multline}
which can be written as
\begin{equation}\label{eqn:lemma:priority:8}
    \eta_1 = -\lambda_{\ell^D(e)}\Bigl(V^{\bm{\eta}}_e(n-1)-V^{\bm{\eta}}_e(n+1)\Bigr) - \lambda_{\ell^D(e)}R_e + \frac{\lambda_{\ell^D(e)}}{\lambda_{\ell^R(e)}}(g_e(\bm{\eta})+\eta_2).
\end{equation}
Together with \eqref{eqn:bellman_2}, we obetain
\begin{equation}\label{eqn:lemma:priority:9}
    \eta_1 = -\lambda_{\ell^D(e)}\Bigl(V^{\bm{\eta}}_e(n) - V^{\bm{\eta}}_e(n+1)\Bigr).
\end{equation}
Based on the definition, for ${\color{black}n}\in\mathscr{N}_+$, if $a^*_1(n) = 0$ then $a^*_1(n+1) = 0$. 
Then, based on \eqref{eqn:bellman_2}, \eqref{eqn:lemma:priority:9} leads to \eqref{eqn:lemma:priority}.

For ${\color{black}n=-N}\in\mathscr{N}_N$ and  $\eta_1=\Gamma(\eta_2)$, we obtain 
$V^{\bm{\eta}}_e({\color{black}-N}) - V^{\bm{\eta}}_e({\color{black}-(N-1)}) = \frac{\lambda_{\ell^D(e)}R_e-g_e(\bm{\eta})-\eta_1}{\lambda_{\ell^D(e)}}$ and $g_e(\bm{\eta}) = 0$.
That is, 
\begin{equation}\label{eqn:lemma:priority:9_5}
    \eta_1 = \lambda_{\ell^D(e)}\Bigl(V^{\bm{\eta}}_e({\color{black}-(N-1)})-V^{\bm{\eta}}_e({\color{black}-N})\Bigr) + \lambda_{\ell^D(e)}R_e.
\end{equation}
For this critical $\eta_1 = \Gamma(\eta_2)$, since $a^*_1({\color{black}-N}) = 0$, $a^*_1({\color{black}-(N-1)})=0$ because of the definition.
It follows that $V^{\bm{\eta}}_e({\color{black}-(N-1)})-V^{\bm{\eta}}_e({\color{black}-N}) = -\frac{\eta_2}{\lambda_{\ell^R(e)}}$. Plugging it in \eqref{eqn:lemma:priority:9_5}, we obtain \eqref{eqn:lemma:priority}.

For ${\color{black}n=0} \in \mathscr{N}_0$, $\eta_1=\Gamma(\eta_2)$ satisfies 
\begin{equation}\label{eqn:lemma:priority:10}
    \eta_1 = -\lambda_{\ell^D(e)} \Bigl(H^{\bm{\eta}}_e({\color{black}0}) - V^{\bm{\eta}}_e{\color{black}1}\Bigr) + \lambda_{\ell^D(e)}R_e,
\end{equation}
and, since $a^*_1({\color{black}0}) = 0 $ leads to $a^*_1({\color{black}1})=0$, 
\begin{equation}\label{eqn:lemma:priority:11}
    V^{\bm{\eta}}_e({\color{black}1}) = \frac{\lambda_{\ell^R(e)}R_e - g_e(\bm{\eta}) - \eta_2}{\lambda_{\ell^R(e)}} + H^{\bm{\eta}}_e({\color{black}0}).
\end{equation}
Plugging \eqref{eqn:lemma:priority:11} in \eqref{eqn:lemma:priority:10}, 
we obtain
\begin{equation}
    \eta_1 = -\frac{\lambda_{\ell^D(e)}}{\lambda_{\ell^R(e)}}(g_e(\bm{\eta}) + \eta_2) + 2\lambda_{\ell^D(e)}R_e \geq -\frac{\lambda_{\ell^D(e)}}{\lambda_{\ell^R(e)}}(g_e(\bm{\eta}) + \eta_2) + \lambda_{\ell^D(e)}R_e.
\end{equation}
It proves the lemma.

\proof{Proof of Proposition~\ref{prop:priority}.}
For $e\in{\color{black}\mE}$, sufficiently small $\eta_2$ and $\eta_1 = -\frac{\lambda_{\ell^D(e)}}{\lambda_{\ell^R(e)}}(g_e(\bm{\eta})+\eta_2) + \lambda_{\ell^D(e)}R_e$, from Lemma~\ref{lemma:priority}, it is optimal to set $a^*_1({\color{black}n}) = 0$ for all ${\color{black}n}\in\tilde{\mathscr{N}}\backslash \{{\color{black}N}\}$.
For such $\bm{\eta}$, since $a^*_1({\color{black}n})=0$ for all $\tilde{\mathscr{N}}\backslash\{{\color{black}N}\}$, $g_e(\bm{\eta}) = 0$.
It follows that, for all ${\color{black}n}\in\tilde{\mathscr{N}}\backslash\{{\color{black}N}\}$, there exists $\Gamma\in\mathscr{B}_1({\color{black}n},e)$ such that
\begin{equation}\label{eqn:prop:priority:1}
\Gamma(\eta_2) = -\frac{\lambda_{\ell^D(e)}}{\lambda_{\ell^R(e)}}\eta_2 + \lambda_{\ell^D(e)}R_e.
\end{equation}
When $\lambda_{\ell^D(e)} = \lambda_{\ell^R(e)}$, \eqref{eqn:prop:priority:1} becomes $\Gamma(\eta_2) = -\eta_2 + \lambda_{\ell^D(e)}R_e$, for which the second item is independent from $\eta_2$.
It proves the proposition.
\Halmos

\section{Simulation Settings}\label{app:uniform_network_settings}

In Section~\ref{sec:example}, we specify the eligibility for a driver and a rider to match each other: 
if both can save costs through a shared ride compared to the costs of the non-sharing rides, the driver and the rider are eligible to match each other; otherwise, not eligible.
In particular, let ${d}_{1}$ and ${d}_{2}$ represent the expected travel distances for the driver and the rider, respectively, and $d^S_{1,2}$ and $d^D_{1,2}$ represent the expected travel distance between the sources and the destinations, respectively, of the two travelers. The rider must drop off first, and let $d=d^S_{1,2}+d_2+d^D_{1,2}$ represent the expected total travel distance for the two travelers if they share a ride. 
The ride charge of going from one location to another is specified to be linear to the travel distance with a slope coefficient $b>0$ (measured in dollars per unit distance) if a traveler does not share a ride with others; and the charge of a ride per unit distance becomes $\gamma b$, $\gamma\in(1,2)$, if two travelers share the ride and share the charge of the entire journey.
When the sum of the expected charge for both travelers going alone is larger than the expected total distance of going together, that is, $b(d_{1}+d_{2}) > \gamma b d$, the two travelers can match each other.
More precisely, if $\ell_1,\ell_2\in{\color{black}\mL}$ represent the labels of the types of the above-mentioned driver and the rider, respectively, then, if $b(d_{1}+d_{2}) > \gamma b d$, there exists a match $(\ell_1,\ell_2)\in\mE$.
In this section, we choose $b=3$ and $\gamma = 1.5$.

For each $e\in{\color{black}\mE}$, the expected reward of the match $R_e$ is set to be $\gamma b (d_{\ell^D(e)}+d_{\ell^R(e)} - d_e)$ where $d_{\ell}$ ($\ell\in{\color{black}\mL}$) is the expected travel distance for the traveler type $\ell$ and $d_e \coloneqq d^S_{\ell^D(e),\ell^R(e)} + d_{\ell^R(e)}+d^D_{\ell^D(e),\ell^R(e)}$ is the expected total travel distance of the shared ride.

The reneging rate of traveler type $\ell$ waiting for the match $e$, $\mu_{\ell,e}$, is assumed to be negatively correlated to the reward $R_e$ and the expected travel distance $d_{\ell}$. 
In particular, let $\mu_{\ell,e}=\exp(-\upsilon R_e - \beta d_{\ell})$ with correlation coefficients $\upsilon,\beta\in\mathbb{R}_+$, which indicates the negative correlation between $\mu_{\ell,e}$ and $d_{\ell}$. Plugging in $R_e$ as a function of $d_{\ell}$, we consider $\beta=0.0189$ and $\upsilon=0.0054$.
Specify the reneging penalty $C_{\ell,e} = -\zeta\ln \mu_{\ell,e}$ with $\zeta \in \mathbb{R}_0$ for $\ell\in{\color{black}\mL}$ and $e\in{\color{black}\mE}$. 
The penalty is then negatively affected by the increasing reneging rate, and is practical for modeling real situations. Travelers who finally fail to take a shared ride and suffer a longer waiting time are likely to be more upset.
As described in Section~\ref{sec:model}, we implement a finite maximal number of travelers waiting for each match $e\in{\color{black}\mE}$, which avoids the explosion of the counting process when the arrival rates of the drivers and riders differ too much. For our simulations presented in this subsection, we set the maximal number of the waiting travelers $N=5$.

\section{Settings for the Melbourne Ride-Sharing Network}\label{app:melbourne_case}

Table~\ref{table:weight_population_melbourne} shows the traveler demands for the Melbourne suburbs whereas Figure~\ref{fig:melbourne_example:demands} depicts the geographic locations of those suburbs.

\begin{table}[!h]
    \caption{Weighted populations for different districts in Melbourne} \label{table:weight_population_melbourne}
    \small
    \begin{tabular}{lr|lr|lr|lr}
\hline
Suburb&\tabincell{c}{Wt.\\ pop.}&Suburb&\tabincell{c}{Wt.\\ pop.}&Suburb&\tabincell{c}{Wt.\\ pop.}&Suburb&\tabincell{c}{Wt.\\ pop.}\\
\hline
Abbotsford & 9848
&
Albert Park & 17197
&
Carlton & 30157
&
Carlton North & 4587
\\
Collingwood & 10372
&
Docklands & 14951
&
East Melbourne & 5739
&
Fitzroy & 11720
\\
Fitzroy North & 13609
&
Kensington & 11862
&
Melbourne & 147291
&
North Melbourne & 25331
\\
Parkville & 8434
&
Port Melbourne & 18129
&
Prahran & 22027
&
Princes Hill & 4587
\\
Richmond & 45039
&
South Melbourne & 17110
&
South Yarra & 41007
&
Southbank & 24376
\\\hline
\end{tabular}
\end{table}









\section{Uniform Network}\label{app:simulation_uninetwork}

\begin{figure}[t]
\centering
\begin{minipage}[]{0.66\textwidth}
\subfigure[]{\includegraphics[width=0.49\linewidth]{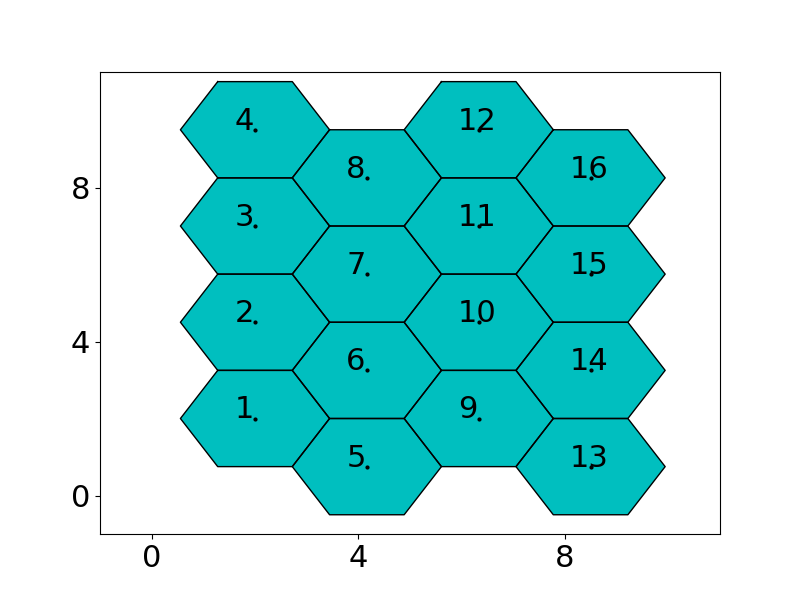}\label{fig:uniform_example:1}}
\subfigure[]{\includegraphics[width=0.49\linewidth]{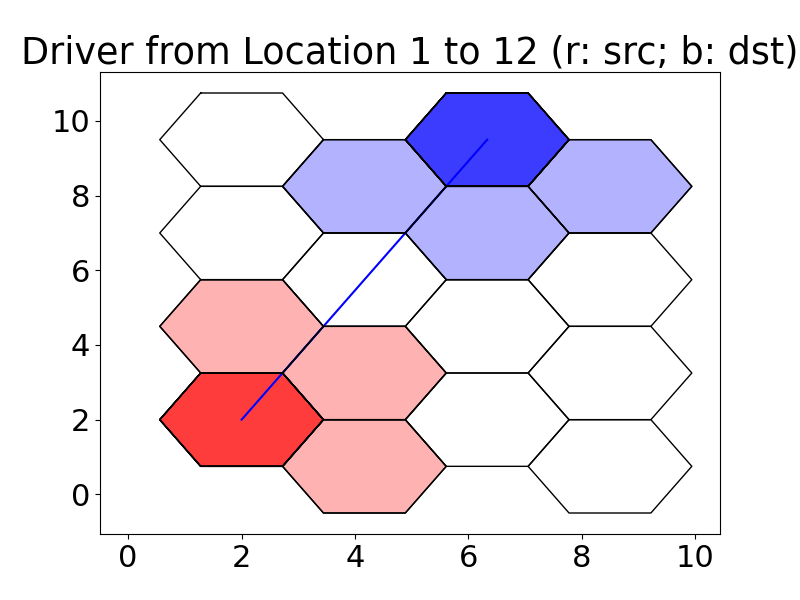}\label{fig:uniform_example:2}}
\vspace{-0.5cm}\caption{Transportation network in the uniform case: (a) all destinations; and (b) eligible traveler types that can match the traveler type $11$.
}\label{fig:uniform_example}
\end{minipage}
\begin{minipage}[]{0.33\textwidth}
\vspace{0.1cm}
\includegraphics[width=\linewidth]{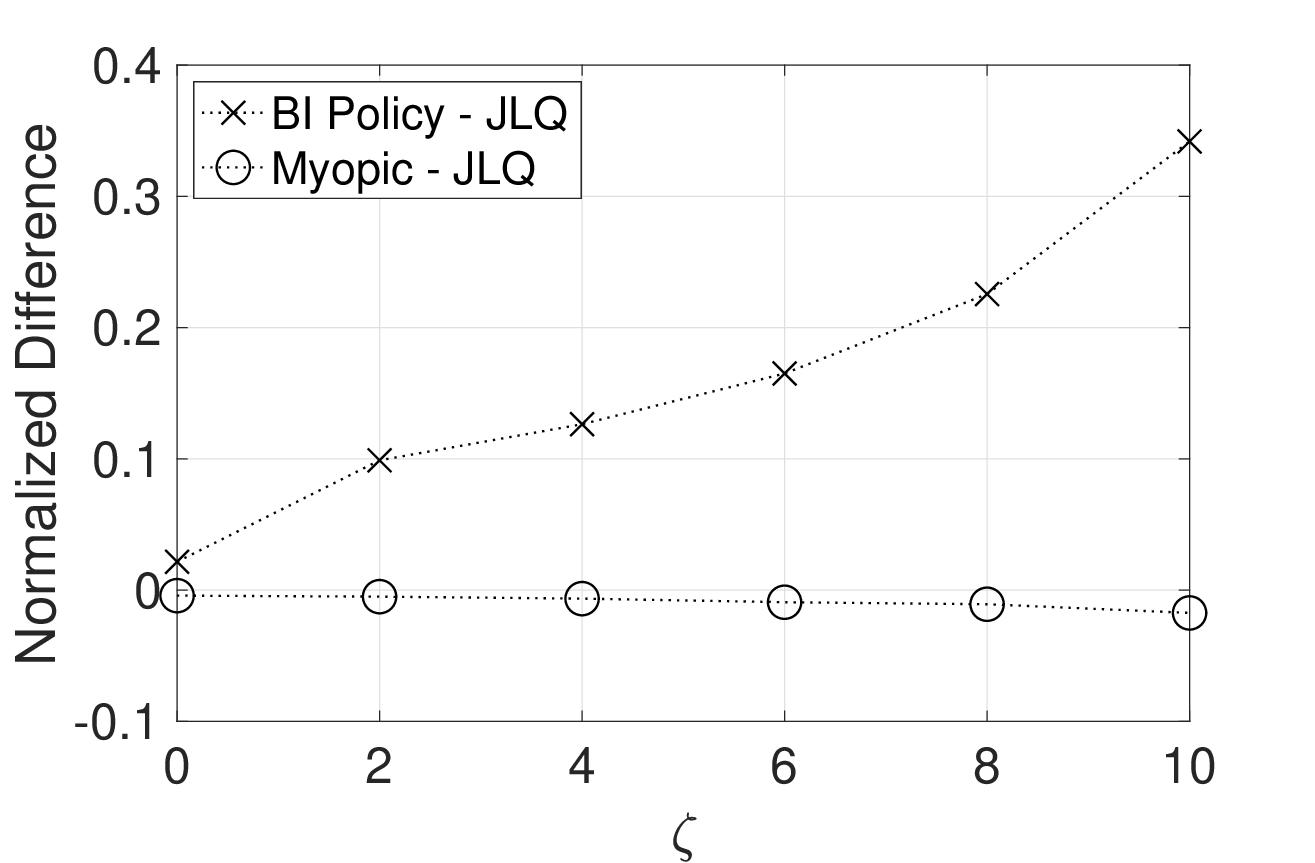}\label{fig:uniform_reward:2}
\caption{The relative difference of specified policy to JLQ with respect to the average reward.}
\label{fig:uniform_reward}
\end{minipage}
\end{figure}

Consider a transportation network with multiple traveler \emph{locations} uniformly deployed in a coordinate.
As illustrated in Figure~\ref{fig:uniform_example:1}, each location is a random variable taking values in a hexagon area, as in Figure~\ref{fig:uniform_example}, of which the travel distance between locations is approximated by the length between their center points.  
The travel distance approximates the expected length of the travel route, which is defined as a random variable in Section~\ref{sec:model}.
As demonstrated in Figure~\ref{fig:uniform_example:1}, here, we consider in total $16$ locations (regions). 
For instance, in Figure~\ref{fig:uniform_example:2}, we illustrate the source and destination locations of a driver type and its available matches.
The red and the blue locations are the source and the destination, respectively, for the corresponding drivers, and the light red and the light blue represent the sources and the destinations, respectively, of the riders that can match the drivers. 
Recall that the drop-off order will affect the availability of matching the corresponding travelers because the rider must drop off first.
In this context, with the $16$ destination areas demonstrated
in Figure~\ref{fig:uniform_example:1}, we consider a case with $240$ types of drivers and $240$ types of riders, where
each driver or rider type corresponds to a pair of source and destination among the $16$ locations.
That is, the number of traveler types $L = 480$.

We consider identical arrival rates of all traveler types, $\lambda_{\ell} = 0.3$ (travelers per minute) for all $\ell\in{\color{black}\mL}$.
The matching mechanism and other system parameters are the same as those for the Melbourne network discussed in Section~\ref{sec:example} (and \ref{app:uniform_network_settings} of the supplementary material), except that, to better match the arrival rates,
the parameters for reneging rates are modified to $\beta = 0.09$ and $\upsilon = 0.03$.
In this context, there are in total $E=682$ matches between the drivers and riders.

Figure~\ref{fig:uniform_reward} demonstrates the relative differences of BI and the myopic policy to JLQ with respect to the average reward as $\zeta$ increases from $0$ to $10$, where the relative differences of BI to JLQ are higher than $10\%$ for all $\zeta\geq 2$, and BI clearly outperforms the myopic policy which achieves similar or even lower average rewards as JLQ.

\begin{figure}
\centering
\begin{minipage}[]{0.33\textwidth}
\centering\vspace{-0.cm}
\includegraphics[width=\linewidth]{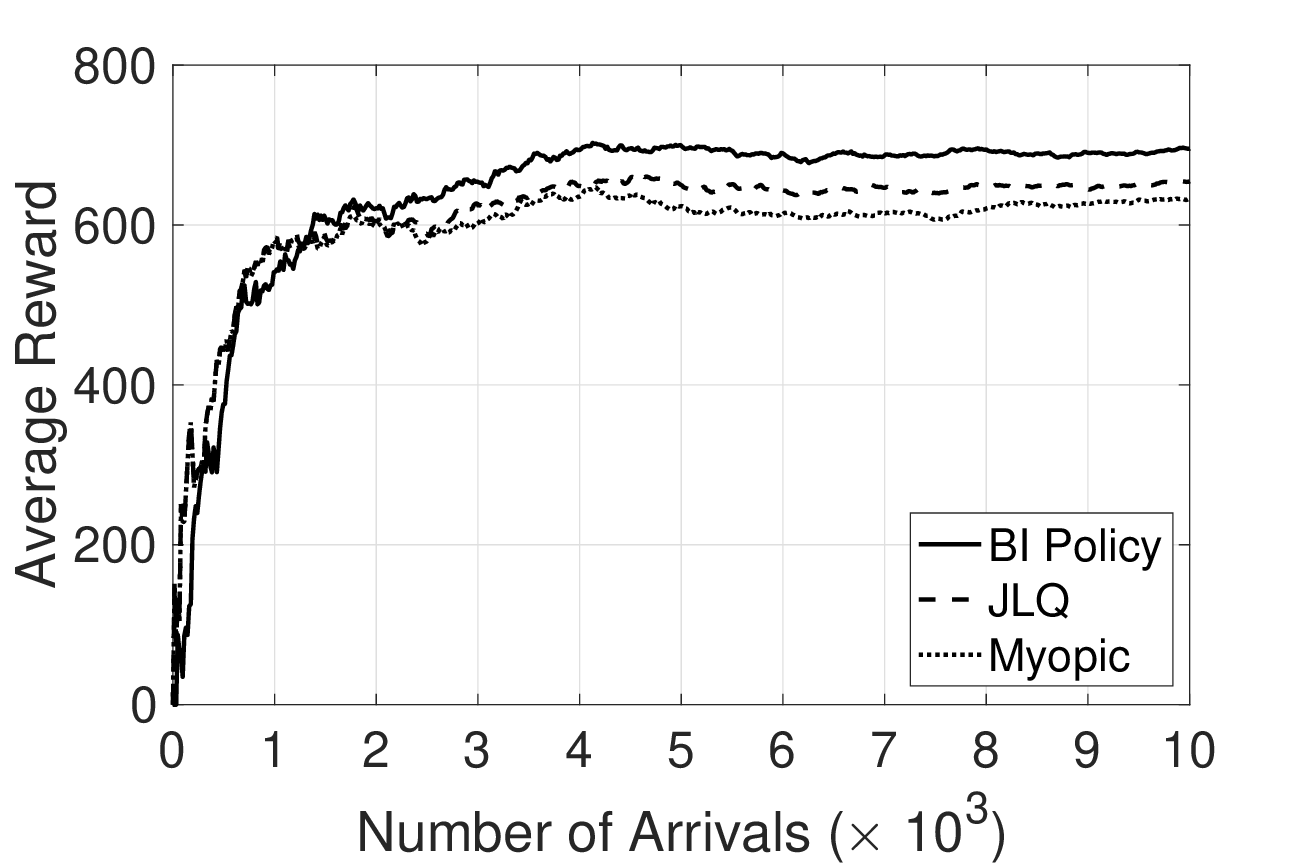}
\caption{ Average reward under different policies against the total number of arrived travelers when $\zeta=4$.}
\label{fig:timeline_arrivals}
\end{minipage}
\begin{minipage}[]{0.66\textwidth}
\centering
\subfigure[]{\includegraphics[width=0.49\linewidth]{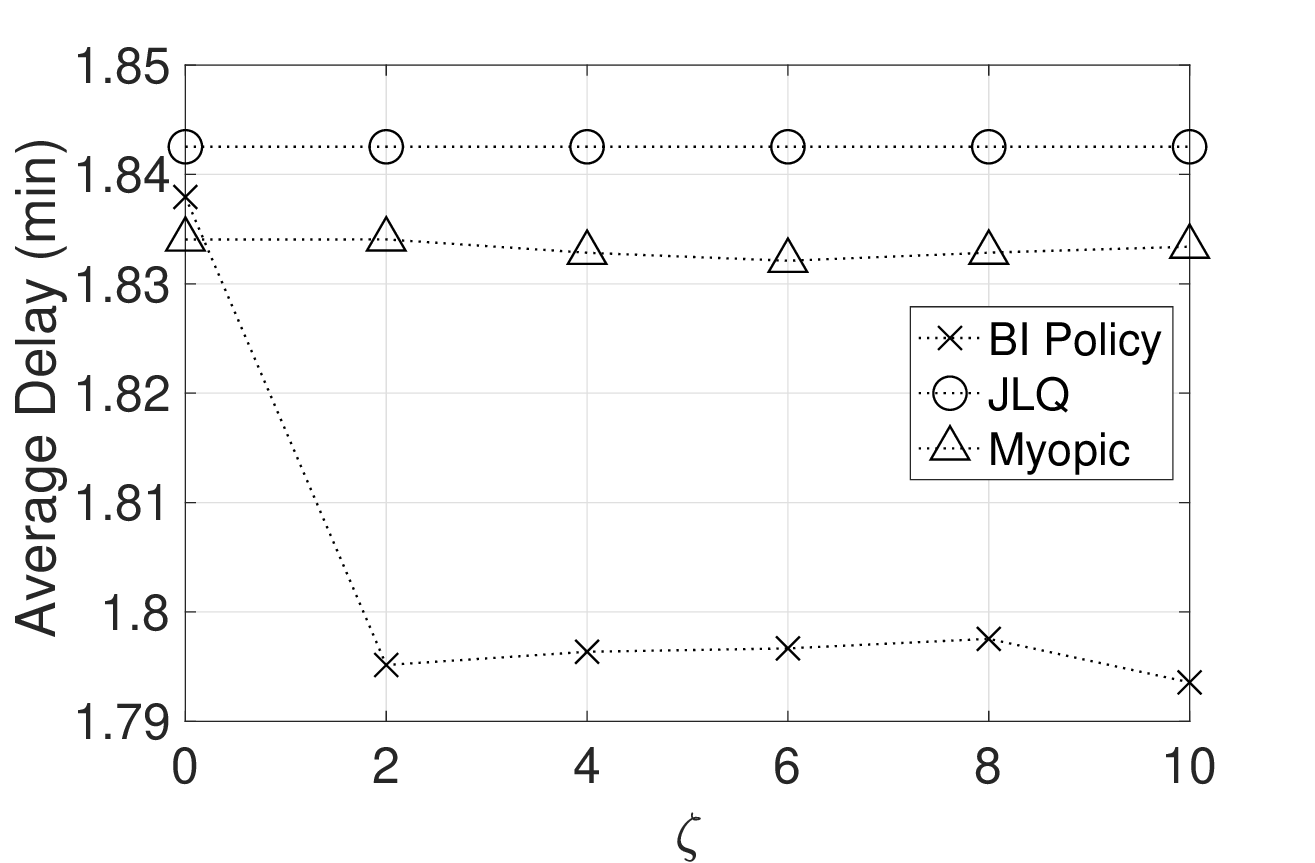}\label{fig:uniform_delay:1}}
\subfigure[]{\includegraphics[width=0.49\linewidth]{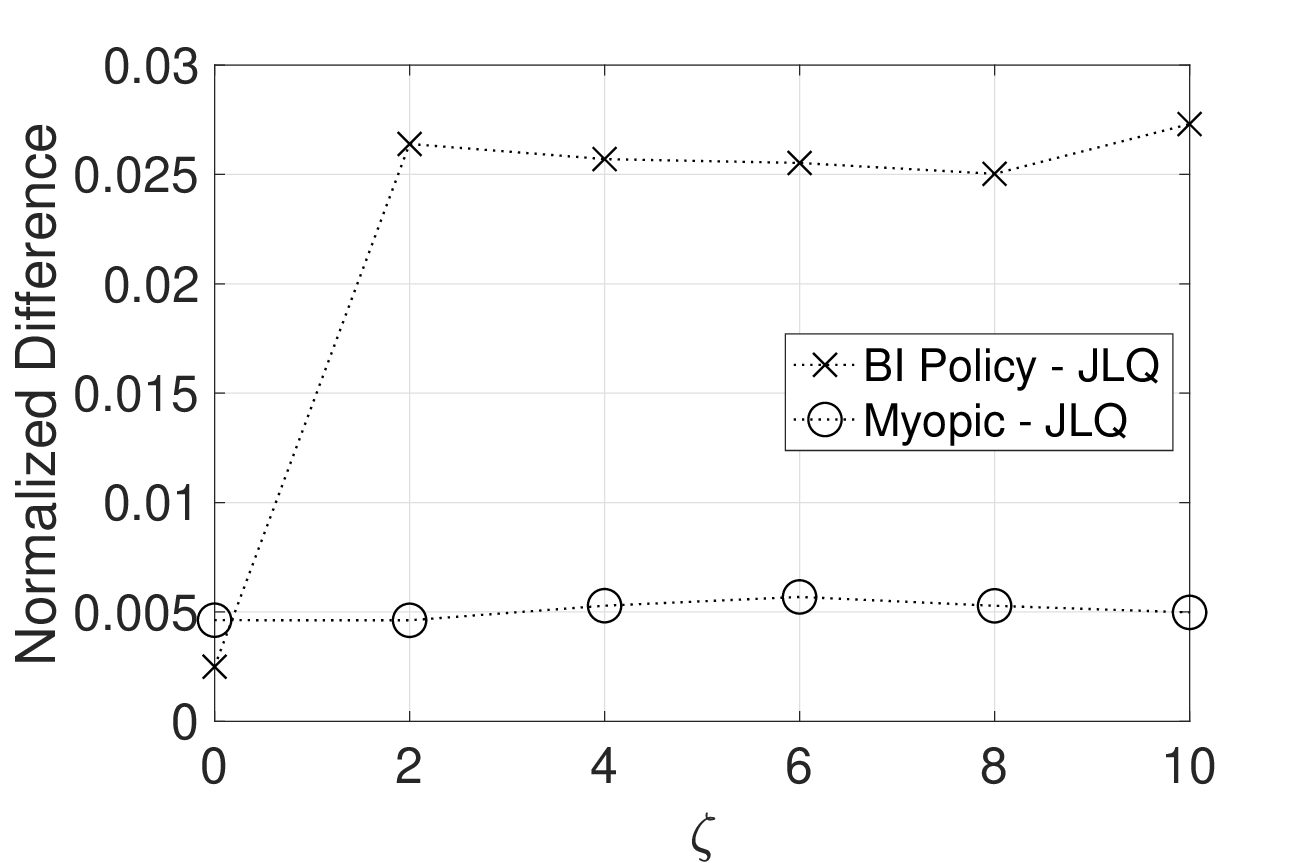}\label{fig:uniform_delay:2}}
\vspace{-0.7cm}\caption{ Average delay of all travelers under different policies against penalty parameter $\zeta$: (a) average delay; and (b) the relative difference of JSQ to the specified policy with respect to the average delay.}
\label{fig:uniform_delay}
\end{minipage}
\end{figure}

Recall that Figure~\ref{fig:uniform_reward} demonstrates the performance in the long-run regime, where the underlying Markov chain of RMP (with finitely many states that communicate with each other) is guaranteed to be stationary. 
For cases with relatively small time horizons, 
we present, in Figure~\ref{fig:timeline_arrivals},  the average rewards of policies against the total number of arrived travelers. The observed curves for all the tested policies become almost flat with a few thousand arrivals. It implies that an urban area with a few thousand ride-sharing travelers per hour needs only one hour to be sufficiently close to the stationary performance.

We plot in Figure~\ref{fig:uniform_delay:1} the average delay of all travelers who once matched someone in our system.
Here, the delay of a traveler means the waiting time of this traveler: the time period lasts from the traveler arrival till he or she either shares a ride or reneges.
Also, define the relative difference of policy $\phi_1$ to policy $\phi_2$ with respect to the average delay as
$\frac{\chi^{\phi_1}-\chi^{\phi_2}}{\chi^{\phi_2}}$,
where $\chi^{\phi}$ represents the average delay of all travelers who once matched with someone under a policy $\phi$.
We plot in Figure~\ref{fig:uniform_delay:2} the relative differences of JLQ to the myopic policy and BI with respect to the average delay. 
In Figure~\ref{fig:uniform_delay}, the average delay is relatively small (less than $2$ minutes) for all the tested policies and the BI policy experiences slightly lower average delays against the myopic policy and JLQ for all tested $\zeta$. That is, given the clear advantages of BI against JLQ and the myopic policy with respect to the average reward, BI achieves a compatible average delay to the other two policies.

\end{document}